\documentclass[11pt]{amsart}
\usepackage[hmargin=3cm,vmargin=3cm]{geometry}

\usepackage[latin1]{inputenc}
\usepackage{xspace,amssymb,amsfonts,euscript}
\usepackage{amsthm,amsmath,amscd,stmaryrd,latexsym}
\usepackage{palatino, euscript}

\usepackage{graphicx}
\usepackage[all]{xy}
\SelectTips{cm}{}

\usepackage{tikz, tikz-cd}
\usetikzlibrary{matrix, arrows, calc}

\usepackage[dvips]{epsfig}
\usepackage{pinlabel}
\usepackage{psfrag}

\newcommand{\ig}[2]{\vcenter{\xy (0,0)*{\includegraphics[scale=#1]{fig/#2}} \endxy}}

\newcommand{\igc}[2]{\begin{center} \includegraphics[scale=#1]{fig/#2} \end{center}}

\IfFileExists{comments.sty}{\usepackage{comments}}

\RequirePackage{color}
\definecolor{myred}{rgb}{0.75,0,0}
\definecolor{mygreen}{rgb}{0,0.5,0}
\definecolor{myblue}{rgb}{0,0,0.65}

\RequirePackage{ifpdf}
\ifpdf
  \IfFileExists{pdfsync.sty}{\RequirePackage{pdfsync}}{}
  \RequirePackage[pdftex,
   colorlinks = true,
   urlcolor = myblue, 
   citecolor = mygreen, 
   linkcolor = myred, 
   pagebackref,
   bookmarksopen=true]{hyperref}
\else
  \RequirePackage[hypertex]{hyperref}
\fi

\newtheorem{thm}{Theorem}[section]
\newtheorem{lemma}[thm]{Lemma}

\newtheorem{prop}[thm]{Proposition}
\newtheorem{cor}[thm]{Corollary}
\newtheorem{claim}[thm]{Claim}  

\newtheorem{goal}[thm]{Goal}

\newtheorem*{prop*}{Proposition}

\theoremstyle{definition}
\newtheorem{defn}[thm]{Definition}
\newtheorem{constr}[thm]{Construction}
\newtheorem{notation}[thm]{Notation}

\newtheorem{example}[thm]{Example}
\newtheorem{exercise}[thm]{Exercise}
\newtheorem{assumption}[thm]{Assumption}

\theoremstyle{remark}
\newtheorem{remark}[thm]{Remark}

\numberwithin{equation}{section}

    \def\BM{{\mathbb{B}}}
\def\CG{{\mathfrak C}}    \def\CM{{\mathbb{C}}}

  \def\gg{{\mathfrak g}}  
\def\HG{{\mathfrak H}}  \def\hg{{\mathfrak h}}

    \def\NM{{\mathbb{N}}}
    
    \def\PM{{\mathbb{P}}}
    
    \def\RM{{\mathbb{R}}}
    \def\SM{{\mathbb{S}}}
    \def\TM{{\mathbb{T}}}

    \def\XM{{\mathbb{X}}}
    
    \def\ZM{{\mathbb{Z}}}

    \def\CC{{\mathcal{C}}}
    \def\DC{{\mathcal{D}}}
    
    \def\FC{{\mathcal{F}}}
    \def\GC{{\mathcal{G}}}
\def\HB{{\mathbf H}}    
    \def\IC{{\mathcal{I}}}
    
    \def\KC{{\mathcal{K}}}

    \def\OC{{\mathcal{O}}}
\def\PB{{\mathbf P}}    \def\PC{{\mathcal{P}}}

    \def\SC{{\mathcal{S}}}

\def\a{\alpha}
\def\b{\beta}

\def\G{\Gamma}
\def\d{\delta}
\def\D{\Delta}

\def\l{\lambda}
\def\L{\Lambda}

\def\o{\omega}
\def\Om{\Omega}
\def\z{\zeta}

\let\phi=\varphi
\let\tilde=\widetilde

\usepackage{bbm}
\def\C{{\mathbbm C}}
\def\N{{\mathbbm N}}
\def\R{{\mathbbm R}}
\def\Z{{\mathbbm Z}}
\def\Q{{\mathbbm Q}}
\def\1{\mathbbm{1}}

\newcommand{\un}{\underline}

\newcommand{\ot}{\otimes}
\newcommand{\pa}{\partial}
\newcommand{\co}{\colon}

\renewcommand{\to}{\rightarrow}
\newcommand{\into}{\hookrightarrow}

\newcommand{\define}{\stackrel{\mbox{\scriptsize{def}}}{=}}

\renewcommand{\sl}{\mathfrak{sl}}

\newcommand{\mt}{\emptyset}

\newcommand{\refequal}[1]{\xy {\ar@{=}^{#1}
(-1,0)*{};(1,0)*{}};
\endxy}

\newcommand{\Kar}{\textbf{Kar}}
\newcommand{\Hom}{{\rm Hom}}

\newcommand{\End}{{\rm End}}

\newcommand{\fin}{{\rm fin}}
\newcommand{\aff}{{\rm aff}}

\newcommand{\Zqq}{\ZM[q^{\pm 1}]}

\newcommand{\SBim}{\SM\textrm{Bim}}
\newcommand{\BSBim}{\BM\SM\textrm{Bim}}
\newcommand{\SSBim}{\SC\SM\textrm{Bim}}
\newcommand{\SBSBim}{\SC\BM\SM\textrm{Bim}}

\newcommand{\mSSBim}{m\SC\SM\textrm{Bim}}
\newcommand{\mSBSBim}{m\SC\BM\SM\textrm{Bim}}

\newcommand{\Rep}{\textrm{Rep}}
\newcommand{\Fund}{\textrm{Fund}}
\newcommand{\RepOm}{\textrm{Rep}^\Om}
\newcommand{\FundOm}{\textrm{Fund}^\Om}
\newcommand{\RepOmq}{\textrm{Rep}^\Om_q}
\newcommand{\FundOmq}{\textrm{Fund}^\Om_q}
\renewcommand{\sc}{\textrm{sc}}
\newcommand{\wt}{\textrm{wt}}
\newcommand{\adj}{\textrm{adj}}
\newcommand{\rt}{\textrm{rt}}
\newcommand{\tG}{\widetilde{\G}}
\newcommand{\rmv}{\Theta}
\newcommand{\trmv}{\tilde{\Theta}}
\newcommand{\Gr}{Gr}
\newcommand{\PPP}{\PB}
\newcommand{\pt}{\star}

\begin{document}

\begin{abstract} We give an interpretation of $\sl_n$-webs as morphisms between certain singular Soergel bimodules. We explain how this is a combinatorial,
algebraic version of the geometric Satake equivalence (in type $A$). We then $q$-deform the construction, giving an equivalence between representations of
$U_q(\sl_n)$ and certain singular Soergel bimodules for a $q$-deformed Cartan matrix.

In this paper, we discuss the general case but prove only the case $n=2,3$. In the sequel we will prove $n \ge 4$. \end{abstract}

\title{Quantum Satake in type $A$: part I}

\author{Ben Elias} \address{Massachusetts Institute of Technology, Boston}

\maketitle


\section{Introduction}
\label{sec-intro}

\subsection{Equivalent equivalences}
\label{subsec-results}

The \emph{geometric Satake equivalence} (or just \emph{geometric Satake} for short) is an equivalence between two symmetric monoidal abelian categories which can be attached to a reductive
algebraic group. In this paper we state a \emph{Soergel Satake equivalence}, an equivalence of (strict) additive $2$-categories associated to a pair of Langlands dual lie
algebras. Furthermore, in type $A$ we state an \emph{algebraic Satake equivalence}, an equivalence between additive $2$-categories living inside the Soergel Satake equivalence. The fact
that these three equivalences imply each other is reasonably straightforward (given the results of Soergel and H\"arterich). The real meat of this paper is computational: an explicit
construction of the algebraic Satake equivalence in type $A$, coming from a presentation of both additive $2$-categories by generators and relations. This gives a new, simple proof of
geometric Satake in type $A$ (but see Remark \ref{rmk:caveat}). Finally, we present a $q$-deformation of the Soergel and algebraic Satake equivalences in type $A$. This
$q$-deformation has no known geometric source at present.

Let $\gg$ be a simple Lie algebra with Langlands dual $\gg^\vee$. Let $\Om$ denote the fundamental group of $\gg^\vee$, a finite abelian group realized as the weight lattice modulo the
root lattice. Representations of $\gg^\vee$ form a semisimple $\Om$-graded-monoidal category, in the sense that each irreducible object has a character (or highest weight) in $\Om$, and
these characters add under taking tensor products and their summands. This additive $\Om$-graded-monoidal category can also be encoded as an additive $2$-category $\RepOm$ with one object
for each element of $\Om$. We give an introduction to these constructions in \S\ref{subsec-om-monoidal}.

Associated to the affine Dynkin diagram $\tG$ of $\gg$, one can construct an additive $2$-category of \emph{singular Soergel bimodules} $\SSBim$ as in \cite{WilSSB}, having one object for
each proper subgraph of $\tG$. Chapter \ref{sec-ssbim} contains an introduction to singular Soergel bimodules. Within $\SSBim$ lies the full sub-$2$-category of \emph{maximally singular
Soergel bimodules} $\mSSBim$, which has objects for each subgraph of $\tG$ isomorphic to the original Dynkin diagram $\G$. Such subgraphs are parametrized naturally by the set $\Om$.

The Soergel Satake equivalence is an equivalence between $\RepOm$ and $\mSSBim$. Actually, $2$-morphisms in $\mSSBim$ are graded vector spaces, so to be more precise, Soergel Satake is a
$2$-functor from $\RepOm$ to $\mSSBim$ which is essentially surjective up to grading shift, faithful, and full onto degree $0$ maps.

The original geometric Satake equivalence (see, for instance, \cite{GinzburgGS}) can be similarly rephrased as an equivalence of $2$-categories between $\RepOm$ and some $2$-category of
perverse sheaves on affine partial flag varieties. Maximally singular Soergel bimodules should be (roughly) thought of as the equivariant hypercohomologies of these perverse sheaves,
though thankfully they have an independent algebraic definition.

Inside $\RepOm$ one can consider the additive sub-$2$-category $\FundOm$ whose $1$-morphisms are tensor products of fundamental representations. There are certain $1$-morphisms in
$\mSSBim$ which correspond to fundamental representations in $\RepOm$, which could be called \emph{fundamental} singular Soergel bimodules. Clearly, $\FundOm$ should be equivalent to the monoidal sub-$2$-category of $\mSSBim$ generated by fundamental singular Soergel bimodules\footnote{More precisely, there should be a functor from $\FundOm$ to this sub-$2$-category which is essentially surjective up to grading shift, faithful, and full onto degree $0$ maps.}, a statement equivalent to Soergel Satake. Sadly, neither $2$-category is well-understood in general.

In type $A$, however, all fundamental representations are miniscule. In this case, the fundamental singular Soergel bimodules are actually easy to describe, being the generating
$1$-morphisms in the $2$-category of \emph{maximally singular Bott-Samelson bimodules} $\mSBSBim$ inside $\mSSBim$. Moreover, in type $A$ both $\FundOm$ and $\SSBim$ have a presentation
by generators and relations using planar diagrams (\cite{CKM} and \cite{EWGR4SSB} respectively), which is why we consider these categories to be ``algebraic." The equivalence between
$\FundOm$ and $\mSBSBim$ in type $A$ will be called the algebraic Satake equivalence, and will be proven by exploiting these presentations.

This paper proves the algebraic Satake equivalence for $\sl_2$ and $\sl_3$. We require several results about singular Soergel bimodules, which for $\sl_2$ are all available in
\cite{EDihedral}, and for $\sl_3$ are proven in the appendix. For $\sl_{n+1}$ with $n \ge 3$, the requisite background will eventually appear in joint work with Williamson \cite{EWGR4SSB},
though we do not know when this manuscript will become available. Moreover, the $n \ge 3$ case also requires a great deal of additional calculation. For these reasons, we postpone the $n
\ge 3$ case to a followup paper.

\begin{remark} \label{rmk:caveat} The feature which is most obscured by the transformation from geometric Satake to algebraic Satake is the symmetric structure on the original monoidal
category. The symmetric structure on perverse sheaves arises in a complicated way from the study of the Beilinson-Drinfeld Grassmannian \cite{MirkVil}, and it is unclear how this is
translated into the language of Soergel bimodules. One could make the argument that geometric Satake without the symmetric structure is a significantly weaker theorem (we discuss this in
more detail in Remark \ref{rmk:caveat2}). However, this argument seems to rely on the Decomposition Theorem \cite{BBD}, while this paper gives a proof of geometric Satake (or rather,
Soergel Satake) without any reliance on the Decomposition Theorem or the related Soergel conjecture. Moreover, making the equivalence explicit on $2$-morphisms is a non-trivial result,
and the ability to discuss the equivalence in a new language (without Tannakian formalism) is useful. In addition, the $q$-deformation below is new. \end{remark}

\begin{remark} \label{rmk:lusztigsupgrade} One can upgrade algebraic Satake to an equivalence of symmetric monoidal categories, using recent work of Lusztig. It is easy to equip
$\FundOm$ with a natural symmetric structure, using the computational formulas (well-known to knot theorists) for the symmetric structure in representation theory. When this paper was
first written, there was no natural way to construct a natural symmetric structure on singular Soergel bimodules. However, Lusztig has since devised such a structure in an update to
\cite{LuszUnequal14}. Technically, his symmetric structure lives on the Soergel bimodules which arise from our singular Soergel bimodules, though it is not hard to adapt the definition to
the singular Soergel bimodules themselves. As far as we can tell, no connection has yet been made between Lusztig's symmetric structure and the Beilinson-Drinfeld Grassmannian. See also
Remark \ref{rmk:caveatquantum} below. \end{remark}

\subsection{The $q$-deformation}
\label{subsec-qdeform}

The $2$-categories $\RepOm$ and $\FundOm$ admit natural $q$-deformations $\RepOmq$ and $\FundOmq$, describing representations of the quantum group $U_q(\gg^\vee)$. Let us restrict
henceforth to type $A$, where $\FundOmq$ also has a known presentation by generators and relations. The surprising fact is that $\SSBim$ also admits a $q$-deformation $\SSBim_q$ compatible
with the Soergel Satake equivalence. We have a functor $\FundOmq \to \mSBSBim_q$ which is fully faithful onto degree $0$ maps, which we call the \emph{quantum algebraic Satake
equivalence}. At the moment, there is no known geometric source for this $q$-deformation of $\SSBim$. Gaitsgory \cite{Gaits} also has a notion of a quantum Satake equivalence, but there is
currently no connection known between his theory and ours.

We now provide a brief description of the $q$-deformation of Soergel bimodules, and a number of remarks. There are some technicalities which we ignore in the introduction;
\S\ref{sec-ssbim} has a more accurate discussion.

Let $W$ be a Coxeter group with simple reflections $S$. In \cite{Soe5} Soergel defines a \emph{reflection faithful representation} of $W$ to be a vector space $\hg$ over a field $\Bbbk$ on
which $W$ acts faithfully, such that an element of $W$ acts by a reflection on $\hg$ (it has one eigenvalue $-1$ and fixes a codimension $1$ hyperplane) if and only if it is a reflection
in $W$. To such a representation, Soergel associates a monoidal category of Soergel bimodules, and Williamson \cite{WilSSB} a $2$-category $\SSBim$ of singular Soergel bimodules. The
lightning definition is this: consider the coordinate ring $R = \textrm{Sym}(\hg^*)$ equipped with its action of $W$. For any finite parabolic subgroup associated to $I \subset S$ one has
a subring $R^I$ of invariant polynomials, and when $I \subset J$ one has $R^J \subset R^I$. Singular Soergel bimodules are defined to be the summands of (grading shifts of) iterated tensor
products of the induction and restriction bimodules between these various rings $R^I$.

In \cite{EWGR4SB} it is explained how to generalize this construction beyond reflection faithful representations. A \emph{realization} of $W$ is a representation $\hg$ which is free over
an arbitrary commutative ring $\Bbbk$, together with a choice of simple coroots $\D^\vee \subset \hg$ and simple roots $\D \subset \hg^*$, satisfying some natural conditions. The pairing
between $\D$ and $\D^\vee$ is encoded in a \emph{(generalized) Cartan matrix} $A$ with entries $(a_{s,t})_{s,t \in S}$ valued in $\Bbbk$. It need not be the case that $\D$ and $\D^\vee$
form bases for their respective spaces, but when they do, the realization is determined by the Cartan matrix. To any realization, \cite{EWGR4SB} and \cite{EWGR4SSB} provide a $2$-category
analogous to $\SSBim$.

To obtain $\SSBim_q$ one begins with a $q$-deformed version of the $\tilde{A}_n$ Cartan matrix, defined over a base ring $\Bbbk = \Z[q,q^{-1}]$, called the \emph{exotic affine $\sl_{n+1}$ Cartan matrix}. For $n \ge 2$ we use
\begin{equation} \label{slnqCartan} \left( \begin{array}{cccccc}
2 & -1 & 0 & \cdots & 0 & -q^{-1} \\
-1 & 2 & -1 & \cdots & 0 & 0 \\
0 & -1 & 2 & \cdots &  & 0 \\
\vdots & \vdots & \vdots & \ddots & -1 & 0 \\
0 & 0 &   & -1 & 2 & -q \\
-q & 0 & 0 & \cdots & -q^{-1} & 2
\end{array} \right). \end{equation}
For $n=1$, we use
\begin{equation} \label{sl2qCartan} \left( \begin{array}{cc}
2 & -(q+q^{-1}) \\
-(q+q^{-1}) & 2 \end{array} \right). \end{equation}
The parameter $q$ can be specialized to a non-zero complex number, so that these exotic matrices also yield one-parameter families of Cartan matrices over $\Bbbk = \CM$. 	These exotic Cartan matrices do not seem to appear in the literature, nor is there yet a satisfactory geometric explanation for them. (However, the corresponding 1-parameter family of representations of the affine Weyl group does appear in the literature, see the discussion in \S\ref{subsec-exoticrealization}). The author came upon them while hunting for a possible
quantum Satake equivalence, inspired by the $\sl_2$ case which is studied in detail in \cite{EDihedral}. An \emph{exotic realization} of affine $\sl_{n+1}$ is a realization of $\tilde{A}_n$ over $\Z[q,q^{-1}]$ or $\CM$ having such a Cartan matrix, and $\SSBim_q$ will be the associated $2$-category.

\begin{remark} Suppose that $\Bbbk$ is a field containing $\RM$. For any two $s \ne t \in S$, let $m_{st}$ denote the order of $st$ in $W$. A \emph{standard} realization of $W$ is one with
a symmetric Cartan matrix valued in $\RM$, satisfying $a_{s,s}=2$ and $a_{s,t} = - \cos(\frac{\pi}{m_{st}})$ for $s \ne t$. When the Coxeter graph of $W$ is a tree, it is easy to argue
that any Cartan matrix is \emph{standardizable}, i.e. conjugate by a diagonal matrix to a standard Cartan matrix. Conjugating by a diagonal matrix corresponds to rescaling $\D$ and
$\D^\vee$, and does not alter the $2$-category of singular Soergel bimodules.

However, when the Coxeter graph of $W$ is not a tree, such as for $\tilde{A}_n$, a realization need not be standardizable. The reader can check that for $n \ge 2$ the exotic Cartan
matrices above are standardizable if and only if $q = \pm 1$. In particular, when $q \in \C^\times \setminus \R^\times$, the complex realization $\hg$ has no real form. It is not
difficult to deduce that the exotic family over $\CM$ exhausts the possible complex Cartan matrices for $\tilde{A}_n$, up to conjugation by diagonal
matrices, from the fact that the Dynkin diagram of $\tilde{A_n}$ has fundamental group isomorphic to $\ZM$.

Because Coxeter graphs of affine Weyl groups in other types are trees, it seems unlikely that one could use a similar construction to produce a quantum algebraic Satake equivalence in
other types. \end{remark}

\begin{remark} The exotic Cartan matrix has determinant $2 - q^2 - q^{-2}$. Therefore, outside of the classical cases $q = \pm 1$, $\D$ and $\D^\vee$ are each linearly independent. Let
$w_1$ denote the reflection corresponding to the longest root in the finite Weyl group (the part with the usual Cartan matrix) and let $s_0$ denote the affine reflection. When $q$ is
specialized to a primitive $2m$-th root of unity, the element $(s_0 w_1)$ has finite order $m$ on the span of $\D^\vee$. Thus the action of $W$ on the span of $\D^\vee$ is not faithful.
Whenever a Cartan matrix is nondegenerate, $\hg$ can be decomposed as a direct sum of the span of $\D^\vee$ and a trivial representation, so that no exotic realization is faithful at a
non-real root of unity. We will discuss this remark further in section \ref{subsec-exoticrealization}. \end{remark}

\begin{remark} Soergel conjectured that for any reflection faithful representation $\hg$ over a field $\Bbbk$ of characteristic zero, the indecomposable Soergel bimodules descend in the
Grothendieck group to the Kazhdan-Lusztig basis of (in this case) the affine Hecke algebra. Williamson has shown \cite{WilSSB} that when this is true, the classes of the indecomposable
singular Soergel bimodules descend to the Kazhdan-Lusztig basis of the affine Hecke algebroid. Soergel's conjecture is proven in \cite{EWHodge} for the standard reflection representation,
but the arguments also apply to other realizations possessing a real form, and satisfying a positivity property.

Now consider the exotic realization when $n \ge 1$, and set $q \ne \pm 1$ to be a root of unity. Quantum algebraic Satake relates maximally singular Soergel bimodules to the
representation theory of $U_q(\sl_{n+1})$, which is no longer semisimple. The lack of semisimplicity prevents the indecomposable bimodules from having the same character they have
generically. Therefore indecomposable bimodules do not descend to the Kazhdan-Lusztig basis. This does not contradict the results of \cite{EWHodge}, as the realization has no real form.
It also does not contradict the Soergel conjecture, as no exotic realization is reflection faithful at a root of unity. \end{remark}

\begin{remark} \label{rmk:caveatquantum} We continue the discussion from Remark \ref{rmk:lusztigsupgrade} above. One might ask whether there is a $q$-deformation of Lusztig's symmetric
structure that produces a natural braiding on maximally singular Soergel bimodules, agreeing under the quantum algebraic Satake equivalence with the standard braiding on
$U_q(\sl_n)$-representations. We consider this a very interesting question (the answer did not appear obvious)! \end{remark}

\subsection{Addendum: further study of the $q$-deformation}
\label{subsec-further}

As mentioned above, the exotic realization of affine $\sl_{n+1}$ is something of a mystery. There is currently no known geometric or representation-theoretic explanation for this
$q$-deformed Cartan matrix; it was fabricated by the author solely to make quantum algebraic Satake work. However, since this paper originally appeared, the exotic realization has taken
on importance in other ways, and we would like to advocate that it is a fundamental object well worth studying, above and beyond the results of this paper. We briefly discuss two recent
applications.

The geometric Satake equivalence is a categorical explanation for a numerical concurrence first observed by Lusztig: that certain multiplicities attached to the Kazhdan-Lusztig basis of
(a singular version of) the Hecke algebra in affine type agreed with multiplicities in the category of representations of the Langlands dual Lie group. Geometric Satake solves this
mystery by declaring that the categories which govern these multiplicities are in fact equivalent.

A similar numerical concurrence is related to the parabolic Kazhdan-Lusztig basis of the affine Hecke algebra, rather than the singular Kazhdan-Lusztig basis. The multiplicities here are
related to the rational representation theory of the corresponding algebraic group in finite characteristic, at least for sufficiently large characteristic. This relationship goes under
many names: Lusztig's conjecture, Andersen's conjecture, James' conjecture. The categorical equivalence which would explain this numerical concurrence has appeared (outside of type $A$,
only conjecturally) in recent work of Riche-Williamson \cite{RicWil}, which contains an excellent introduction to the topic and many additional references. These same multiplicities are
also related to the representation theory of quantum groups at roots of unity. In this case, the numerical equality has been a theorem for some time, due to the work of many: Kazhdan,
Lusztig, Kashiwara, Tanisaki, Soergel, etcetera. A quantum Riche-Williamson conjecture would give an equivalence of categories between tilting modules in the trivial block of the quantum
group $U_q(\sl_{n+1})$ at a root of unity, and a parabolic version of the diagrammatic category of Soergel bimodules coming from the exotic realization of affine $\sl_{n+1}$. For $\sl_2$
this was proven already by Andersen-Tubbenhauer \cite{AndTub}, to which we refer the reader for a further introduction to this topic. The general case is under investigation.

We were careful to refer to the ``diagrammatic category of Soergel bimodules'' in the previous paragraph, instead of just the category of Soergel bimodules. This diagrammatic category
was constructed in \cite{EWGR4SB}, and agrees with Soergel's algebraic definition when Soergel bimodules behave well (e..g for reflection faithful representations). For arbitrary
realizations, it is the diagrammatic category that is the appropriate generalization. The reflection representation of affine $\sl_{n+1}$ in finite characteristic, and its exotic
realization at a root of unity, are certainly not faithful, so that the diagrammatic category is the category whose multplicities agree with those in the affine Hecke algebra, not
Soergel's category. However, it is still well worth asking: what does the algebraic category of Soergel bimodules look like? What does it categorify? When $q$ is specialized to a $2m$-th
root of unity, the exotic realization factors through the finite complex reflection group $G(m,m,n)$, and Soergel bimodules in this setting appear to be related to this complex
reflection group.

An early investigation of this story is currently underway between the author and Benjamin Young. We study the \emph{exotic NilCoxeter algebra}, the algebra generated by Demazure operators
when $q$ is specialized to a $2m$-th root of unity. Unlike the case of generic $q$, this is a finite dimensional algebra, which we can present by generators and relations (at least for
small $n$, so far). It appears to be an entirely new algebra, with an unfamiliar Poincar\'e polynomial and even a surprising dimension! For example, the exotic NilCoxeter algebra for
$G(2,2,3)$ is 36 dimensional, even though $G(2,2,3) \cong S_4$ only has size 24. We expect these exotic NilCoxeter algebras to be extremely interesting. They have desirable properties: a
unique longest element, for instance. This longest element can be used to define a Frobenius extension structure between the polynomial ring of the exotic realization and the subring of
polynomials invariant under the action of $G(m,m,n)$. This Frobenius extension structure, in turn, is essential to the study of Soergel bimodules.

\subsection{Organization of the paper}
\label{subsec-organization}

In \S\ref{sec-sl2-diag} we define a $2$-category $\FundOm(\sl_2)$ by generators and relations, using the language of planar diagrammatics. We similarly define a graded $2$-category
$\mSBSBim(\sl_2)$. Then we define a $2$-functor $\FundOm(\sl_2) \to \mSBSBim(\sl_2)$, and prove that it is essentially surjective (up to grading shift) and fully faithful (onto degree zero
$2$-morphisms). In \S\ref{sec-sl3-diag}, we repeat the same process for $\sl_3$. These two sections constitute the real mathematical content of this paper (which is quite easy). They are
presented in a vacuum, as it were, without any reference to Soergel bimodules, perverse sheaves, or representation theory, and thus are accessible without any background.

In the subsequent chapters we fill in the details, eventually connecting this computational result to geometric Satake. In \S\ref{sec-webs} we give an introduction to the diagrammatic
approach to $\Rep(\gg^\vee)$. We also describe the group $\Om$ and the notion of $\Om$-graded-monoidal categories and $\Om$-$2$-categories. In \S\ref{sec-ssbim} we give background on
singular Soergel bimodules, and state the Soergel Satake equivalence. Between these two chapters, both sides of the algebraic Satake equivalence described in \S\ref{sec-sl2-diag} and
\S\ref{sec-sl3-diag} will be explained.

Finally, in \S\ref{sec-reformulate} we explain how to connect Soergel Satake to geometric Satake. This takes place in three steps, with two being very straightforward, and the last relying
upon the difficult work of Soergel and H\"arterich.

With the exception of \S\ref{sec-reformulate}, everything is already written with the $q$-deformation built in. This leads to a host of complications when defining Soergel bimodules. As a
result, \S\ref{sec-ssbim} is not the easiest introduction to Soergel bimodules in general, though it does provide an in-depth introduction to the subtleties involved when dealing with
``odd-unbalanced" realizations, such as the exotic affine realization of $\sl_n$.

This paper is organized with two audiences in mind: the neophyte and the expert. For the neophyte with little to no experience with geometric Satake, we suggest reading the paper in order.
We hope that the concrete, combinatorial approach will make it readable, and the hands-on experience of \S\ref{sec-sl2-diag} and \S\ref{sec-sl3-diag} will help when trying to understand
the more abstract approaches. The paper should be entirely accessible until \S\ref{sec-reformulate}. We do assume the reader is familiar with the technology of planar diagrammatics for
$2$-categories with adjunction: an introduction can be found in \cite[section 2]{LauDiagrams}. The expert who is more interested in the connection from geometric Satake to Soergel Satake
and then to algebraic Satake is welcome to skip directly to \S\ref{sec-reformulate}, only backtracking when indicated for several definitions.

\textbf{Acknowledgements}

The author is extremely grateful to Geordie Williamson, George Lusztig, Pavel Etingof, Roman Bezrukavnikov, Joel Kamnitzer, and Shrawan Kumar for a variety of useful discussions.

\section{The diagrammatic approach: $\sl_2$}
\label{sec-sl2-diag}

Most of the material in this chapter can be found in \cite{EDihedral}. We include it as a warm-up for the more complicated case of $\sl_3$, and $\sl_n$ in the sequel.

In this chapter, let $\Om = \{ \overline{0},\overline{1}\}$ denote the group $\Z / 2\Z$. Let the set $S=\{b,r\}$ be an $\Om$-torsor, and assign the elements colors: blue to $b$ and red to
$r$.

We also assign names and colors to proper subsets $I \subsetneq S$. We abusively let $r$ (resp. $b$) denote the singleton subset $\{r\} \subsetneq S$. Let $\mt$ denote the empty set, to
which we assign the color white.

\begin{notation} \label{not:present} In this chapter and the next, we will define several $2$-categories with similar presentations. The objects will form a finite set, and each object
will have an assigned color. There will be a set of generating $1$-morphisms, with the property that for any two colors $s_1, s_2$ there is at most one generating $1$-morphism from $s_1$
to $s_2$. Therefore, a 1-morphism can be represented uniquely by the (non-empty) sequence $\un{x} = s_1 s_2 \ldots s_m$ of colors through which it passes (though not all sequences of
colors are permitted). We read 1-morphisms from right to left, so that $\un{x}$ has source $s_m$ and target $s_1$. The identity $1$-morphism of an object $s$ is thus also denoted $s$. We
represent a composition of 1-morphisms diagrammatically as a sequence of dots on the line, separating regions of different colors. \end{notation}

\begin{example} \label{example:1mor2} The $1$-morphism $brbrbr$ in the $2$-category $\FundOmq$ below: $\ig{1}{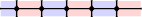}$. \end{example}

For the rest of this chapter we fix a $\Zqq$-algebra $\Bbbk$. In order to connect these results to representation theory we will set $\Bbbk = \Q(q)$, or $\Bbbk = \C$ with $q=1$. However,
the algebraic Satake equivalence is defined in more generality.

\subsection{Webs}
\label{subsec-webs-sl2}

Background on this material can be found in \cite{GooWen, WestburyTL}. More specifics on this particular 2-colored version can be found in \cite[section 2]{ELib} or \cite[section 4]{EDihedral}.

\begin{defn} \label{defn:webs-sl2} Let $\FundOmq = \FundOmq(\sl_2)$ (also known as the $2$-colored Temperley-Lieb category) be the $\Bbbk$-linear $2$-category defined as follows, using
Notation \ref{not:present}. It has objects $S= \{r,b\}$, and has generating $1$-morphisms $rb$ and $br$. Thus Example \ref{example:1mor2} gives a $1$-morphism in $\FundOmq$.

The 2-morphisms are generated by colored cups and caps. More precisely, there is a cap map $brb \to b$ and a cup map $b \to brb$, as pictured below, as well as
the corresponding maps with the colors switched. \igc{1}{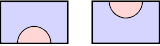}

There are two types of relations, which also hold with the colors switched. Recall that $[2]=q+q^{-1}$.

The \textbf{Isotopy relation}:
\begin{equation} \label{eq:isotopy} \ig{1}{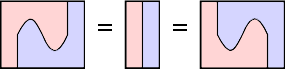} \end{equation}
	
The \textbf{Circle relation}:
\begin{equation} \label{eq:circle} 	{
	\labellist
	\small\hair 2pt
	 \pinlabel {$-[2]$} [ ] at 60 20
	\endlabellist
	\centering
	\ig{1}{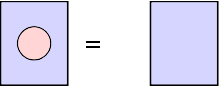}
	} \end{equation}
This ends the definition.
\end{defn}

The existence of cups and caps, together with \eqref{eq:isotopy}, is equivalent to the statement that the 1-morphisms $rb$ and $br$ are biadjoint. There is an action of $\Om$ on $\FundOm$
which permutes the colors.

\begin{defn} A \emph{crossingless matching} is an isotopy class of colored 1-manifold with boundary embedded in the planar strip, without any closed components, providing a matching of
boundary points. \end{defn}

\begin{example} A crossingless matching in $\Hom(brbrbrbrbrb,brbrbrbrb)$: \igc{1}{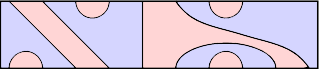} \end{example}
	
\begin{claim} Each morphism space $\Hom(\un{x},\un{y})$ has a basis over $\Bbbk$ given by crossingless matchings. \end{claim}

Now we define another $2$-category.

\begin{defn} Let $\RepOmq = \RepOmq(\sl_2)$ denote the following full sub-$2$-category of $\textbf{Cat}$. The objects are in bijection with $\Om$, where $\overline{0}$ is identified with
the category of even $U_q(\sl_2)$ representations, and $\overline{1}$ with odd $U_q(\sl_2)$ representations. The $1$-morphisms are (functors arising from) tensor products with finite
dimensional $U_q(\sl_2)$ representations, and the $2$-morphisms are all natural transformations (i.e. $U_q(\sl_2)$-morphisms between tensor products). \end{defn}

Now we fix an identification of $\Om$ with $S$. There is a $2$-functor from $\FundOmq$ to $\RepOmq$ which on objects sends $S$ to $\Om$. It sends both $1$-morphisms $rb$ or $br$ to the
tensor product with the standard representation $V$. The cups and caps are sent to inclusions and projections respectively between $V \ot V$ and the trivial summand $\Lambda^2 V$. These
inclusions and projections are unique up to scalar; one chooses the scalars so that \eqref{eq:isotopy} holds. For any such choice of scalars, it turns out that \eqref{eq:circle} will also
hold.

\begin{claim} Suppose that $\Bbbk = \Q(q)$. For any identification of $S$ with $\Om$, this $2$-functor is well-defined and fully faithful. \end{claim}

The proofs of these claims are well-known, and are reasonable exercises for the uninitiated. The result also holds when $\Bbbk = \CM$, $q = 1$, and $U_q(\sl_2)$ is replaced by $\sl_2$. A
similar claim holds when $\Bbbk = \Z[q,q^{-1}]$, for the correct integral form of $U_q(\sl_2)$.

\subsection{Singular Soergel diagrams}
\label{subsec-ssbimdiag-sl2}

For another introduction to this material, see \cite[section 5]{EDihedral}. We now make a mild assumption on $\Bbbk$: that the ideal generated by $2$ and $[2]$ is the unit ideal (c.f.
Demazure Surjectivity in \cite[section 3.3]{EDihedral}).

\begin{defn} \label{defn:SBSBim-sl2} Let $\SBSBim_q = \SBSBim_q(\sl_2)$ be the graded $\Bbbk$-linear $2$-category defined as follows, using Notation \ref{not:present}. The objects are proper subsets $I
\subsetneq S$, and the generating $1$-morphisms are $\{ r \mt, \mt r, b \mt, \mt b\}$.

\begin{example} The $1$-morphism $b \mt r \mt r \mt b \mt$: $\ig{1.5}{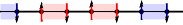}$. \end{example}

In the drawing above, we placed an orientation on the dots consistent with the rule that the right-hand color contains the left-hand color. This orientation is redundant information, but
helps to make $2$-morphisms easier to describe. Similarly, we will color the strands in our pictures by the color which is added or subtracted by the $1$-morphism.

The $2$-morphisms are generated by oriented cups and caps (again, the orientation itself is redundant). There are cap maps $\mt b \mt \to \mt$ and $b \mt b \to b$, as well as cup maps in
the other direction, as pictured below. There are also maps with blue and red switched. \igc{1}{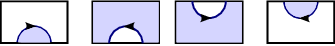}

We place a grading on the $2$-morphisms, where clockwise cups and caps have degree $+1$, and anticlockwise cups and caps have degree $-1$. The relations between $2$-morphisms are listed
below (with the color-switched versions assumed). They are all homogeneous.

The \textbf{Isotopy relation}:
\begin{equation} \label{eq:isotopysbsb} \ig{1}{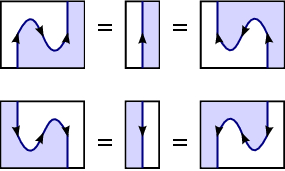} \end{equation}

The \textbf{Empty Circle relation}:
\begin{equation} \label{ccwcirc} \ig{1}{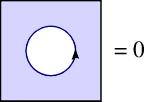} \end{equation}

The \textbf{Cartan relations}:
\begin{subequations} \label{thecartanrelations}
\begin{equation} \label{circeval} \ig{1}{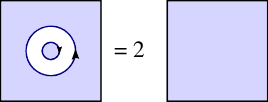} \end{equation}
\begin{equation} \label{circeval1} 	{
	\labellist
	\small\hair 2pt
	 \pinlabel {$-[2]$} [ ] at 75 26
	\endlabellist
	\centering
	\ig{1}{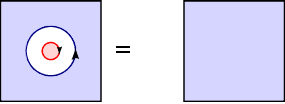}
	}
 \end{equation}
\end{subequations}

The \textbf{Forcing relations}:
\begin{subequations} \label{thecircleforcingrelations}
\begin{equation} \label{circforcesame} \ig{1}{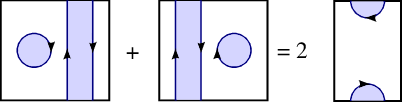} \end{equation}
\begin{equation} \label{circforce1} 	{
	\labellist
	\small\hair 2pt
	 \pinlabel {$-[2]$} [ ] at 75 25
	 \pinlabel {$[2]$} [ ] at 143 26
	\endlabellist
	\centering
	\ig{1}{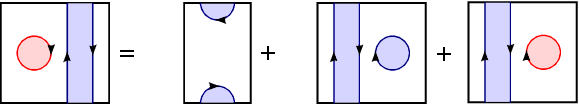}
	}
 \end{equation}
\end{subequations}
This ends the definition.
\end{defn}

The existence of cups and caps, together with \eqref{eq:isotopysbsb}, is equivalent to the statement that the 1-morphisms $b\mt$ and $\mt b$ are biadjoint. There is an action of $\Om$ on
$\FundOm$ which permutes the colors; however, this is special to the $\sl_2$ case.

\subsection{The equivalence}
\label{subsec-equiv-sl2}

\begin{defn} Let $\mSBSBim_q$ denote the full sub-$2$-category of $\SBSBim_q$ whose objects are $\{r,b\}$. In other words, we allow $1$-morphisms and $2$-morphisms which contain the color
white, but their right-hand and left-hand colors must be either red or blue. \end{defn}

\begin{defn} Let $\FC$ be the $2$-functor $\FundOmq \to \mSBSBim_q$ defined as follows. On objects, it sends red to red and blue to blue. On $1$-morphisms, it sends $rb \mapsto r \mt b$
and $br \mapsto b \mt r$. On $2$-morphisms, it acts as below. Visually, the map on 2-morphisms takes a crossingless matching and widens each strand into a white region, with its boundary
oriented anticlockwise. \igc{1}{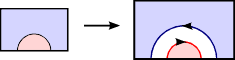} \end{defn}

\begin{claim} The $2$-functor $\FC$ is well-defined, and its image consists of degree $0$ maps. \end{claim}

\begin{proof} The isotopy relation in $\FundOmq$ follows from the isotopy relations in $\SBSBim_q$. The circle relation in $\FundOmq$ follows from the Cartan relations in $\SBSBim_q$.
\end{proof}

\begin{defn} A $2$-functor $\GC$ from a $\Bbbk$-linear $2$-category $\CC$ to a $\Bbbk$-linear graded $2$-category $\DC$ is a \emph{degree-zero equivalence} if the following properties
hold. \begin{itemize} \item $\GC$ induces a bijection between the objects. \item $\GC$ induces an isomorphism $\Hom_{\CC}(\un{x},\un{y}) \to \Hom_{\DC}^0(\GC(\un{x}),\GC(\un{y}))$ to the
morphisms of degree $0$. \item $\Hom_{\DC}^k(\GC(\un{x}),\GC(\un{y}))=0$ for all $k<0$. \item Every $1$-morphism in $\DC$ is isomorphic to a direct sum of grading shifts of $\GC(x)$
for some $1$-morphism $x$ in the Karoubi envelope $\Kar(\CC)$; i.e. $\GC$ is \emph{essentially surjective up to grading shift}. \end{itemize} \end{defn}

\begin{thm} \label{mainthm-sl2} Suppose that the ideal generated by $2$ and $[2]$ in $\Bbbk$ is the unit ideal. Then the $2$-functor $\FC$ is a degree-zero equivalence. \end{thm}

This is one of the main results of \cite[section 4]{EDihedral}, where its proof is dispersed throughout many pages. We will sketch the proof after some more material.

\subsection{Singular Soergel diagrams and polynomials}
\label{subsec-ssbimdiag-sl2-alt}

Let us introduce some important notation for certain $2$-morphisms in $\SBSBim_q$. Let $\a_b$ (resp. $\a_r$) denote a clockwise blue (resp. red) circle. These are endomorphisms of the
identity $1$-morphism $\mt$. Therefore we have a homomorphism from the polynomial ring $R = \Bbbk[\a_b,\a_r]$ to $\End(\mt)$, which is a graded homomorphism provided we set the degree of
$\a_b$ to be $+2$. We use this homomorphism to replace disjoint unions of circles by ``boxes" (labeled with polynomials).

{
\labellist
\small\hair 2pt
 \pinlabel {$\a_b^2 \a_r$} [ ] at 84 20
\endlabellist
\centering
\igc{1}{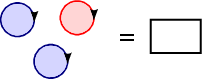}
}

We define endomorphisms $b$ and $r$ of $R$ by the formulas \begin{subequations} \begin{eqnarray} b(\a_b) = -\a_b, & \qquad & b(\a_r) = \a_r + [2] \a_b\\ r(\a_b) = \a_b + [2] \a_r, & \qquad
& r(\a_r) = -\a_r. \end{eqnarray} \end{subequations} In particular, $b^2 = r^2 = 1$. This action on linear terms is a $q$-deformation of the reflection representation of the affine Weyl
group of $\sl_2$.

We can define subrings $R^b$ and $R^r$ of invariant polynomials. In \cite[section 5.2.1]{EDihedral} it is shown how to define a map from $R^b$ to the endomorphism ring of the identity
$1$-morphism $b$. Thus we can place a box labeled by a polynomial within a blue region as well, provided the polynomial is invariant under the operator $b$.

One crucial fact about these rings is that $R^b \subset R$ is a \emph{Frobenius extension}. Roughly, this means that $R$ is free as an $R^b$-module, and that the functors of induction and
restriction between $R$-modules and $R^b$-modules are biadjoint. Part of the data of a Frobenius extension is an $R^b$-linear map $\pa_b \co R \to R^b$, in this case defined by the formula
\[ \pa_b(f) = \frac{f - b(f)}{\a_b}.\] Moreover, $R$ has dual bases (as an $R^b$-module) $\{a_i\}$ and $\{a_i^*\}$ with respect to $\pa_b$, in the sense that \[\pa_b(a_i a_j^*) = \d_{ij}.
\] (These dual bases are not part of the data of a Frobenius extension.) There is a canonical ``coproduct" element $\D_b = \sum_i a_i \ot a_i^*$ living in $R \ot_{R^b} R$, independent of
the choice of dual bases. We often use Sweedler notation $\D_b = \D_{b (1)} \ot \D_{b (2)}$. One invariant of a Frobenius extension is the product $\mu_b = \Delta_{(1)} \Delta_{(2)} =
\sum_i a_i a_i^* \in R$, which in this case is equal to $\a_b$. We provide more detail on Frobenius extensions in
\S\ref{subsec-frobenius}.

The values of $\pa_s(\a_t)$ for $s,t \in S$ are encoded in a \emph{Cartan matrix}, yielding precisely the $q$-deformed affine Cartan matrix of $\sl_2$, given in \eqref{sl2qCartan}. These values are also evident from the Cartan relation \eqref{thecartanrelations}.

Note that $\pa_b$ is surjective, because it is $R^b$-linear and $\pa_b(a_i a_i^*)=1$. This relied on our assumption that the ideal in $\Bbbk$ generated by $2$ and $[2]$ contains the unit.
After all, $\pa_b(\a_b) = 2$ and $\pa_b(\a_r) = -[2]$, so that this ideal comprises all the scalars in the image of $\pa_b$.

Now we give an alternate description of $\SBSBim_q$.

\begin{defn} \label{defn:SBSBim-sl2-alt} Let $\SBSBim_q$ denote the $2$-category with objects and $1$-morphisms as in Definition \ref{defn:SBSBim-sl2}. The $2$-morphisms will be generated
by cups and caps, as well as boxes. A box appearing in a region labeled $I \subset S$ is labeled by a polynomial in $R^I$, and has degree equal to the degree of that polynomial (with the
convention $\deg \a_b = 2$). Boxes add and multiply as polynomials do. In addition to the isotopy relation \eqref{eq:isotopysbsb} we have the following relations.

\begin{subequations} \label{singularrelations}
\begin{equation} \label{circleis} 	{
	\labellist
	\small\hair 2pt
	 \pinlabel {$\a_b$} [ ] at 70 18
	 \pinlabel {$\a_r$} [ ] at 178 18
	\endlabellist
	\centering
	\ig{1}{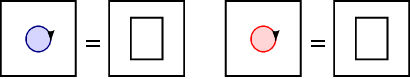}
	}. \end{equation}
\begin{equation} \label{slidepolyintoblue} 	{
	\labellist
	\small\hair 2pt
	 \pinlabel {$f$} [ ] at 18 18
	 \pinlabel {$f$} [ ] at 147 18
	\endlabellist
	\centering
	\ig{1}{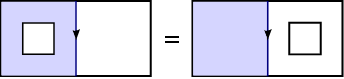}
	} \textrm{ when } f \in R^b. \end{equation}
\begin{equation} \label{demazureis} 	{
	\labellist
	\small\hair 2pt
	 \pinlabel {$f$} [ ] at 20 20
	 \pinlabel {$\pa_b(f)$} [ ] at 83 20
	\endlabellist
	\centering
	\ig{1}{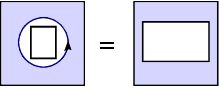}
	}. \end{equation} \begin{equation} \label{circforcesamealt} {
	\labellist
	\small\hair 2pt
	 \pinlabel {$\D_b$} [ ] at 85 23
	\endlabellist
	\centering
	\ig{1}{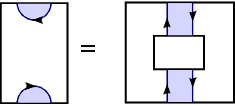}
	}. \end{equation}
\end{subequations}
In this last equation, $\D_b$ is meant to represent the action of $R \ot_{R^b} R$ inside $\End(\mt b \mt)$, by placing $\D_{(1)}$ in the white region on the left, and $\D_{(2)}$ in the white region on the right. Equation \eqref{slidepolyintoblue} guarantees that this action is well-defined. \end{defn}

In \cite{EDihedral} it is shown that this definition of $\SBSBim_q$ is equivalent to the one in Definition \ref{defn:SBSBim-sl2}. It is also a reasonable exercise for the reader.

Now we define the singular Soergel $2$-category.

\begin{defn} Let $\SSBim_q = \SSBim_q(\sl_2)$ denote the following full sub-$2$-category of $\textbf{Cat}$. The objects are proper subsets of $S$, where $I \subsetneq S$ is identified with
the category of graded $R^I$-modules (using the notation $R^\mt = R$). The $1$-morphisms are summands of (grading shifts of) iterated induction and restriction functors, between $R^b$- or
$R^r$-modules and $R$-modules. The $2$-morphisms are all natural transformations (i.e. homogeneous bimodule maps). \end{defn}

There is a $2$-functor from $\SBSBim_q$ to $\SSBim_q$ which is the identity on objects. It sends $b \mt$ to the restriction functor (with a grading shift by $1$), and $\mt b$ to the
induction functor. The cups and caps are sent to units and counits of adjunction.

\begin{claim} Suppose that $\Bbbk = \Q(q)$, or that $q=1$ and $\Bbbk$ is a field of characteristic $\ne 2$. Then this $2$-functor is well-defined and fully faithful. \end{claim}

See \cite[section 5.2.2 and Corollary 5.30]{EDihedral} for more details.

\subsection{Sketch of proof}
\label{subsec-sketch-sl2}

Let us quickly sketch the proof of Theorem \ref{mainthm-sl2} given in \cite[section 4]{EDihedral}. There are three main facets of the argument. One uses idempotent decompositions (or the
algebraic theory of Frobenius extensions) to prove essential surjectivity up to grading shift. One uses considerations on the Grothendieck group to determine the graded dimension of morphism
spaces. Finally, one uses a diagrammatic argument to show that $\FC$ is full.

The fact that $R$ is free over $R^b$ says that the $R^b$-bimodule $R$ should split into copies of the $R^b$-bimodule $R^b$. In $\SBSBim_q$, this amounts to the fact that $b \mt b$ is
isomorphic to several copies of $b$ (with grading shifts). One way to prove this splitting is to construct an idempotent decomposition of the identity of $b \mt b$. This is obtained by
rotating \eqref{circforcesamealt} by 90 degrees. When $2$ is invertible, this can be realized more concretely by rotating \eqref{circforcesame} by 90 degrees and dividing by $2$.

Any $1$-morphism in $\mSBSBim_q$ which alternates between red and blue, i.e. $r \mt b \mt r \mt \ldots \mt b$, is clearly in the image of $\FC$. Any other $1$-morphism must contain either
$b \mt b$ or $r \mt r$. By the above paragraph, one can replace $b \mt b$ with $b$ after taking summands, and similarly with $r \mt r$. Therefore, any object is a direct sum of alternating
$1$-morphisms, proving the essential surjectivity up to grading shift.

Using standard representation theory or combinatorics, one knows the dimension of the space of crossingless matchings between any two colored sequences $\un{x}$ and $\un{y}$ in $\FundOmq$.
Using the fact that $\SSBim_q$ categorifies the Hecke algebroid with its standard trace, one knows the (graded) dimension of morphism spaces between any two objects in $\mSBSBim_q$. From
this one can deduce that there are no morphisms of negative degree between alternating $1$-morphisms, and that the maps in degree zero are the correct size for $\FC$ to induce
isomorphisms. (Note that the proof in \cite{EDihedral} does not use precisely the same argument: it uses Theorem \ref{mainthm-sl2} to prove that $\SSBim_q$ categorifies the Hecke algebra,
not the other way around. Instead, one proves computationally that $\FC$ is faithful.)

Thus we need only show that $\FC$ is full onto maps of minimal degree. This is a purely diagrammatic argument, the details of which are not particularly relevant. First one proves that any
closed diagram (i.e. a diagram in the endomorphism ring of an identity $1$-morphism) reduces to a polynomial (in $R$ or $R^b$ or $R^r$, depending on the color of the boundary). Obviously
any diagram with a non-trivial polynomial is not of minimal degree, so we can assume our diagram has no closed subdiagrams. Now we wish to show that blue regions are only adjacent
(crossing white regions) to red regions, which would imply that one could deformation retract the white regions into a crossingless matching, and the diagram is in the image of $\FC$.
However, if there are two separate blue regions separated by a white region, then ``fusing" them \igc{1}{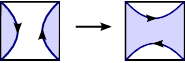} would yield a morphism of lower degree. (This argument is analogous to the
one made in \cite{EDihedral}, which passes through a separate diagrammatic calculus to achieve the same result.)

\section{The diagrammatic approach: $\sl_3$}
\label{sec-sl3-diag}

In this chapter, let $\Om = \{ \overline{0},\overline{1},\overline{2}\}$ denote the group $\Z / 3\Z$. Let the set $S=\{b,r,y\}$ be an $\Om$-torsor, and assign the elements colors: blue to
$b$, red to $r$, and yellow to $y$. The action of $+1 \in \Om$ acts in alphabetical order to send $b$ to $r$, $r$ to $y$, and $y$ to $b$. We call this \emph{color rotation}.

We also assign names and colors to proper subsets $I \subsetneq S$. We abusively let $r$ (resp. $b$, $g$) denote the singleton subset $\{r\} \subsetneq S$. Let $\mt$ denote the empty set,
to which we assign the color white. To pairs we associate the natural compound color: purple $p$ to $\{r,b\}$, green $g$ to $\{b,y\}$, and orange $o$ to $\{r,y\}$. Therefore, the action of
$+1 \in \Om$ sends $p$ to $o$, $o$ to $g$, and $g$ to $p$ (reverse alphabetical order, unfortunately).

For the rest of this chapter we fix a $\Zqq$-algebra $\Bbbk$. In order to connect these results to representation theory we will set $\Bbbk = \Q(q)$, or $\Bbbk = \C$ with $q=1$. However,
the algebraic Satake equivalence is defined in more generality.

\subsection{Webs}
\label{subsec-webs-sl3}

Background on this material can be found in \cite{Kup} or \cite{CKM}.

\begin{defn} \label{defn:webs-sl3} Let $\FundOmq = \FundOmq(\sl_3)$ (also known as the $2$-category of colored $\sl_3$ webs) be the $\Bbbk$-linear $2$-category defined as follows, using
Notation \ref{not:present}. It has objects $\{o,g,p\}$ (i.e. $I \subset S$ of size $2$), and has generating $1$-morphisms from any color to any different color.

\begin{example} The $1$-morphism $pogogp$: $\ig{1}{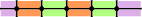}$. \end{example}

The 2-morphisms are generated by colored cups and caps, and by trivalent vertices. More precisely, there is a cap map $sts \to s$ and a cup map $s \to sts$ for any two colors $s,t$. There
is a trivalent vertex $ogpo \to r$, and another trivalent vertex $opgo \to r$, as pictured below. \igc{1}{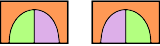} The relations are as follows; they hold for every valid coloring. Recall that $[2] = q+q^{-1}$
and $[3] = q^2 + 1 + q^{-2}$.

The \textbf{Isotopy relations}:
\begin{equation} \label{eq:isotopy3} \ig{1}{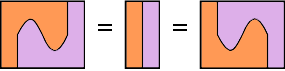} \end{equation}

\begin{equation} \label{eq:cyclicity} \ig{1}{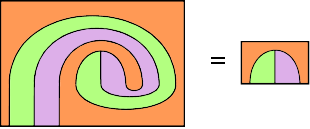} \end{equation}

The \textbf{Circle relation}:
\begin{equation} \label{eq:circle3} 	{
	\labellist
	\small\hair 2pt
	 \pinlabel {$[3]$} [ ] at 60 20
	\endlabellist
	\centering
	\ig{1}{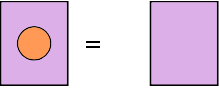}
	} \end{equation}
	
The \textbf{Bigon relation}:
\begin{equation} \label{eq:bigon} 	{
	\labellist
	\small\hair 2pt
	 \pinlabel {$-[2]$} [ ] at 60 20
	\endlabellist
	\centering
	\ig{1}{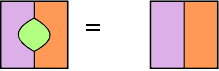}
	} \end{equation}
	
The \textbf{Square relation}:
\begin{equation} \label{eq:square} \ig{1}{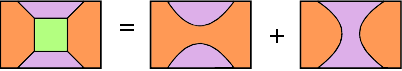} \end{equation}
This ends the definition. \end{defn}

The existence of cups and caps, together with \eqref{eq:isotopy3}, is equivalent to the statement that the 1-morphisms $po$ and $op$ are biadjoint. The cyclicity relation
\eqref{eq:cyclicity} states that all $2$-morphisms are cyclic with respect to these biadjunctions. Together, this allows one to unambiguously define the various trivalent vertices
appearing above as rotations of the generating trivalent vertices. There is an action of $\Om$ on $\FundOm$ by color rotation.

\begin{defn} A \emph{non-elliptic web} is an isotopy class of colored trivalent graph with boundary embedded in the planar strip, which may have interior hexagons, octagons, etcetera, but
may have no interior squares or bigons. \end{defn}

\begin{example} A non-elliptic web in $\Hom(opgpgpopogo,opgpogopgo)$: \igc{1}{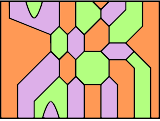} \end{example}

\begin{claim} \label{nonelliptic} Each morphism space $\Hom(\un{x},\un{y})$ has a basis over $\Bbbk$ given by non-elliptic webs. \end{claim}
	
Now we define another $2$-category. Note that there is a bijection between $\Om$ and the set of central characters of $SL_3$.

\begin{defn} Let $\RepOmq = \RepOmq(\sl_3)$ denote the following full sub-$2$-category of $\textbf{Cat}$. The objects are in bijection with $\Om$, where each $\xi \in \Om$ is identified
with the category of $U_q(\sl_3)$ representations having central character $\Om$. The $1$-morphisms are (functors arising from) tensor products with finite dimensional $U_q(\sl_3)$
representations, and the $2$-morphisms are all natural transformations (i.e. $U_q(\sl_3)$-morphisms between tensor products). \end{defn}

Now we fix an identification of $\Om$ with $\{o,g,p\}$. Say we identify $o$ with the trivial character, $g$ with the character of $V_{\o_1}$, and $p$ with the character of $V_{\o_2}$. There
is a $2$-functor from $\FundOmq$ to $\RepOmq$ which on objects sends $\{o,g,p\}$ to $\Om$. The $1$-morphisms $go$, $pg$, and $op$ all correspond to tensoring with $V = V_{\o_1}$, while
their biadjoints correspond to tensoring with $V^* = V_{\o_2}$. The cups and caps correspond to inclusion and projection between $V \ot V^*$ and the trivial summand. The trivalent vertex
$opgo \to o$ corresponds to projection from $V^{\ot 3}$ to the trivial summand $\Lambda^3 V$. These projections and inclusions are unique up to scalar; one chooses the scalars so
that \eqref{eq:isotopy} and \eqref{eq:bigon} hold. For any such choice of scalars, it turns out that \eqref{eq:circle3} and \eqref{eq:square} also hold.

\begin{claim} Suppose that $\Bbbk = \Q(q)$. For any identification of $S$ with $\Om$, this $2$-functor is well-defined and fully faithful. \end{claim}

One again, the proofs of these claims are well-known (also see \cite{Kup}), and are worthwhile exercises. The result also holds when $\Bbbk = \CM$, $q = 1$, and $U_q(\sl_3)$ is replaced by
$\sl_3$. A similar claim holds when $\Bbbk = \Z[q,q^{-1}]$, for the correct integral form of $U_q(\sl_3)$.

\subsection{Polynomials}
\label{subsec-polynomials-sl3}

Our next goal is to define a $2$-category $\SBSBim_q(\sl_3)$, in similar fashion to $\SBSBim_q(\sl_2)$. However, because polynomials are slightly more complex for $\sl_3$ than $\sl_2$, we
choose to discuss them first in this section. Again, we need a mild assumption on $\Bbbk$ in order for the inclusion of an invariant subring to be a Frobenius extension. In this case, the
assumption that $3$ is invertible will suffice.

Let $R = \Bbbk[\a_r, \a_b, \a_y]$ be a polynomial ring in 3 variables, graded so that linear terms have degree $+2$. We define an action of three operators $b,r,y$ on $R$ as follows:
\begin{subequations} \begin{eqnarray} b(\a_b) = -\a_b, & b(\a_r) = \a_r + \a_b, & b(\a_y) = a_y + q^{-1} \a_b,\\ r(\a_b) = \a_b + \a_r, & r(\a_r) = -\a_r, & r(\a_y) = \a_y + q \a_r, \\
y(\a_b) = \a_b + q \a_y, & y(\a_r) = \a_r + q^{-1} \a_y, & y(\a_y) = -\a_y. \end{eqnarray} \end{subequations} In particular, each operator is an involution, and the reader can check that
$(br)^3 = (by)^3 = (ry)^3 = 1$. This is a $q$-deformation of the reflection representation of the affine Weyl group of $\sl_3$. Note that, unlike the case of $\sl_2$, there is no color
symmetry in this action, so the colors are not interchangeable.

Let $I \subsetneq S$, and let $R^I$ denote the subring of polynomials invariant under all reflections in $I$. For $s \in S$ we let $\pa_s \co R \to R^s$ be defined as before, via the formula \begin{equation} \label{eq:partial} \pa_s(f) = \frac{f - s(f)}{\a_s}. \end{equation} This entirely describes the action of the reflection on the linear terms. The values of $\pa_s(\a_t)$ for $s,t \in S$ can be encoded in the $q$-deformed affine $\sl_3$ Cartan matrix.
\begin{equation} \label{sl3qCartan} \left( \begin{array}{ccc}
2 & -1 & -q^{-1} \\
-1 & 2 & -q \\
-q & -q^{-1} & 2
\end{array} \right). \end{equation}

The reader should confirm that $\pa_s$ and $\pa_t$ do not satisfy the braid relation! In particular, one has \begin{subequations} \label{demazurebraid} \begin{eqnarray} \pa_b \pa_r \pa_b &
= & \pa_r \pa_b \pa_r, \\ q \pa_b \pa_y \pa_b & = & \pa_y \pa_b \pa_y, \\ q^{-1} \pa_r \pa_y \pa_r & = & \pa_y \pa_r \pa_y. \end{eqnarray} \end{subequations} Either $\pa_b \pa_y \pa_b$ or
$\pa_y \pa_b \pa_y$ will define a Frobenius structure map $R \to R^g$, but there is no particular reason a priori to choose one over the other. The choice we make will lead to some
convenient formulas.

\begin{defn} \label{defn:roots3} For each $I \subsetneq S$ we will define a set of \emph{positive roots} $\Phi_I$, in such a way that $\Phi_I \subset \Phi_J$ whenever $I \subset J$. We
have $\Phi_{\mt} = \mt$, and $\Phi_s = \{\a_s\}$ for $s \in S$. For doubletons, we have \[ \Phi_{p} = \{\a_b, \a_r, b(\a_r) = \a_r + \a_b\}, \Phi_{g} = \{\a_b, \a_y, b(\a_y) = \a_y +
q^{-1} \a_b\}, \Phi_{o} = \{\a_r, \a_y, r(\a_y) = \a_y + q \a_r\}. \] We call the union of these positive roots the \emph{finite positive roots} of $S$. Let $\mu^I_J$ denote the product of
the roots in $\Phi_J$ which are not in $\Phi_I$.

For later reference, we let $c_{s,t}$ denote the coefficient of $\a_s$ in the third root of $\Phi_{s,t}$. Thus $c_{r,b} = c_{b,r} = c_{y,b} = c_{y,r} = 1$, while $c_{b,y} = q^{-1}$ and
$c_{r,y} = q$. \end{defn}

\begin{defn} \label{defn:frobmaps} For each $I \subsetneq J \subsetneq S$ we will define a map $\pa^I_J \co R^I \to R^J$. These maps are \emph{compatible} in that $\pa^J_K \circ \pa^I_J =
\pa^I_K$ whenever $I \subsetneq J \subsetneq K \subsetneq S$. In particular, we need only define the map $\pa^I_J$ when $I$ and $J$ differ by a single element; this will define all
$\pa^I_K$ uniquely by the comptability requirement, once we check that this is consistently defined. For each $s \in S$ we choose $\pa^\mt_s = \pa_s$ as defined above. Then we
choose \[ \pa^b_{p} = \pa_b \pa_r, \qquad \pa^r_{p} = \pa_r \pa_b, \qquad \pa^b_{g} = q \pa_b \pa_y, \qquad \pa^y_{g} = \pa_y \pa_b, \qquad \pa^r_{o} = q^{-1} \pa_r \pa_y, \qquad \pa^y_{o}
= \pa_y \pa_r. \] The compatibility relations follow from \eqref{demazurebraid}.

Note that $\pa_s \pa_t = c_{s,t} \pa^s_{s,t}$. \end{defn}

\begin{claim} The maps $\pa^I_J \co R^I \to R^J$ equip the ring extension $R^J \subset R^I$ with the structure of a Frobenius extension. For any choice of dual bases, one has $\mu^I_J =
\sum_i a_i a_i^*$. In particular, this implies that $\pa^I_J (\mu^I_J)$ is an integer, equal to the rank of the extension. \end{claim}

This claim is not difficult to prove, by explicitly constructing dual bases. We invite the reader to check that whenever $s \in S$ and $I$ is a doubleton containing $s$, one has
$\pa^s_I(\mu^s_I) = 3$. When performing such computations, one should make use of the \emph{twisted Leibniz rule} \[\pa_s(fg) = \pa_s(f) g + s(f) \pa_s(g). \] Once again, the fact that one
has a Frobenius extension relies upon the fact that $\pa^s_I$ is surjective, which was implied by our assumption that $3$ is invertible.

For more information on choosing structure maps for collections of Frobenius extensions, see sections \S\ref{subsec-frobenius} and \S\ref{subsec-invts}.

Before moving on to the diagrammatic definition of singular Bott-Samelson bimodules, we state some remarkable identities which will play a role in the algebraic Satake equivalence, and
which will generalize in some sense to $\sl_n$.

\begin{claim} The following identities hold. They are defined in a color-symmetric way, but each case should be checked individually. \begin{itemize} \item Let $I,J \subset
S$ be two distinct doubletons, with intersection $s$. Then \begin{equation}\label{howtogetq3} \pa^s_I(\mu^s_J) = [3] = q^2 + 1 + q^{-2}. \end{equation} \item Let $I = \{s,t\}$ be a
doubleton with complement $u \in S$. Let $\mu_I^{s,t}$ denote the product of the roots in $\Phi_I$ which are not in $\Phi_s$ or $\Phi_t$; this is none other than the third root in
$\Phi_I$, either $b(\a_r)$, $b(\a_y)$, or $r(\a_y)$. Then \begin{equation} \label{howtogetq2} \pa_u(\mu_I^{s,t}) = -[2] = -(q+q^{-1}). \end{equation} \item For any $s,t,u \in S$ distinct we have \begin{equation} \label{howtoget1} c_{s,u} c_{t,u} = 1. \end{equation} \end{itemize} \end{claim}

\begin{proof} These are straightforward exercises. \end{proof}

\subsection{Singular Soergel diagrams}
\label{subsec-ssbimdiag-sl3}

We now define $\SBSBim_q(\sl_3)$ analogously to the $\sl_2$ case in Definition \ref{defn:SBSBim-sl2-alt}.

\begin{defn} \label{defn:SBSBim-sl3} Let $\SBSBim_q = \SBSBim_q(\sl_3)$ be the graded $\Bbbk$-linear $2$-category defined as follows, using Notation \ref{not:present}. The objects are
proper subsets $I \subsetneq S$. There is a generating $1$-morphism between two colors if and only if one contains the other and their sizes differ by one.

\begin{example} The $1$-morphism $b \mt r o y g b \mt$: $\ig{1.5}{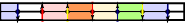}$. \end{example}

In the drawing above, we placed an orientation on the dots consistent with the rule that the right-hand color contains the left-hand color. This orientation is redundant information, but
helps to make $2$-morphisms easier to describe. Similarly, we will color the strands in our pictures by the color which is added or subtracted by the $1$-morphism.

The $2$-morphisms are generated by oriented cups and caps (again, the orientation itself is redundant) and boxes, as well as crossings. Within a region colored $I$ one can place a box
labeled by a polynomial in $R^I$. There are crossing maps $\mt b p \to \mt r p$ as pictured below. More generally, there are crossing maps $\mt s I \to \mt t I$ for any doubleton $I=\{s,t\} \subset S$. \igc{1}{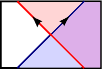}

We place a grading on the $2$-morphisms as follows. A clockwise cup between $\mt$ and $s \in S$ has degree $+1$, while a clockwise cup between $s$ and $\{s,t\}$ has degree $+2$.
Anticlockwise cups have the opposite degree. Upward-oriented crossings have degree $0$. Note however that one can use cups and caps to draw a crossing with the arrows pointing either
down or sideways. By computation, a downward-oriented crossing also has degree $0$, while a sideways-oriented crossing has degree $+1$. Boxes have degree given by the degree of the
polynomial.

The (homogeneous) relations given below are in terms of the Frobenius structures defined above. As a result, they will not be invariant under color rotation! For illustrative purposes,
however, we choose specific colors to exemplify relations.

One has the isotopy relations as well as all the relations of \eqref{singularrelations}, with the obvious modification
\begin{equation} \label{circleis3}
{
\labellist
\small\hair 2pt
 \pinlabel {$\a_b$} [ ] at 70 18
 \pinlabel {$\a_r$} [ ] at 178 18
 \pinlabel {$\a_y$} [ ] at 286 18
\endlabellist
\centering
\ig{1}{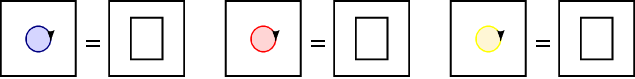}
}
\end{equation}

We also have the analogous relations for the other Frobenius extensions.
\begin{subequations} \label{morefrobenius}
\begin{equation} \label{bigdemazure} 	{
	\labellist
	\small\hair 2pt
	 \pinlabel {$f$} [ ] at 20 21
	 \pinlabel {$\pa^b_p(f)$} [ ] at 89 21
	\endlabellist
	\centering
	\ig{1}{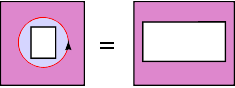}
	} \end{equation}
\begin{equation} \label{bigcircle} 	{
	\labellist
	\small\hair 2pt
	 \pinlabel {$\mu^b_p$} [ ] at 70 20
	\endlabellist
	\centering
	\ig{1}{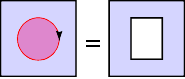}
	} \end{equation}
\begin{equation} \label{slidepolyintopurple} 	{
	\labellist
	\small\hair 2pt
	 \pinlabel {$f$} [ ] at 18 19
	 \pinlabel {$f$} [ ] at 147 19
	\endlabellist
	\centering
	\ig{1}{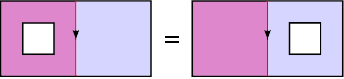}
	} \textrm{ when } f \in R^p \end{equation}
\begin{equation} \label{bigcircforce} 	{
	\labellist
	\small\hair 2pt
	 \pinlabel {$\D^b_p$} [ ] at 86 24
	\endlabellist
	\centering
	\ig{1}{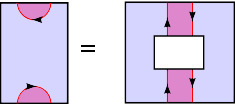}
	} \end{equation}
\end{subequations}

We have a number of relations analogous to the Reidemeister II move.
\begin{subequations} \label{R2moves}
\begin{equation} \label{singR2} \ig{1}{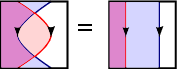} \end{equation}
\begin{equation} \label{singR2nonoriented1} 	{
	\labellist
	\small\hair 2pt
	 \pinlabel {$\pa \D_p$} [ ] at 68 16
	\endlabellist
	\centering
	\ig{1}{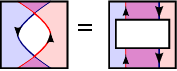}
	} \end{equation}
\begin{equation} \label{singR2nonoriented2} 	{
	\labellist
	\small\hair 2pt
	 \pinlabel {$\mu^{r,b}_p$} [ ] at 71 18
	\endlabellist
	\centering
	\ig{1}{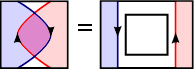}
	} \end{equation}

These relations require some explanation. The symbol $\pa \D_p$ in \eqref{singR2nonoriented1} represents the element \[\pa \D_p = \pa_b(\D^r_{p (1)}) \ot \D^r_{p (2)} = \D^b_{p (1)} \ot \pa_r (\D^b_{p (2)}) \in R^b \ot_{R^p} R^r. \] Meanwhile, the symbol $\mu^{r,b}_p$ is the product of the roots in $\Phi_p$ which are not in either $\Phi_b$ or $\Phi_r$, which in this particular case is the linear term $b(\a_r)$.

We pause to note a consequence of the preceding relations, which is akin to \eqref{singR2nonoriented1}.
\begin{equation} \label{singR2nonoriented1var} {
\labellist
\small\hair 2pt
 \pinlabel {$f$} [ ] at 21 16
 \pinlabel {$\D^b_{p (1)}$} [ ] at 81 16
 \pinlabel {$\pa_r(f \D^b_{p (2)})$} [ ] at 149 16
\endlabellist
\centering
\ig{1}{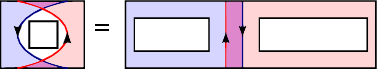}
} 
\end{equation}
\end{subequations}

The relations above are actually completely general, for any compatible system of Frobenius extensions. See \cite{EWFrob} for more details. The final relations below are specific to our
particular circumstance.
\begin{equation} \label{squaresbsb}
	{
	\labellist
	\small\hair 2pt
	 \pinlabel {$q^{-1}$} [ ] at 120 73
	 \pinlabel {$q^{}$} [ ] at 120 22
	\endlabellist
	\centering
	\ig{1}{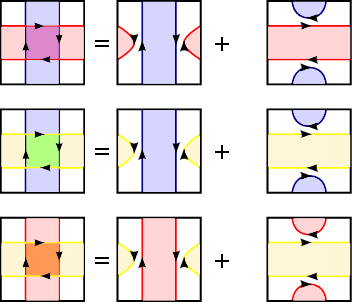}
	}
\end{equation}
To compare these relations to our choice of Frobenius structure, the reader should note the following fact. Take the RHS of each equation, and apply a cap to each of the four walls of the diagram, to obtain diagrams with 3 bubbles each. The cubic polynomial obtained is precisely $\mu_I$ for the doubleton color $I$ appearing on the LHS.  Said another way, the coefficient of the diagram with colors $s,t$ and more connected components colored $s$ is $c_{s,t}$.

This ends the definition. The diagrams which arise in this fashion are called \emph{singular Soergel diagrams}.\end{defn}

There is also a ``boxless" definition analogous to Definition \ref{defn:SBSBim-sl2}, though it is less convenient. In this alternate definition, the 2-morphisms are given purely in terms
of oriented colored $1$-manifolds, and boxes containing polynomials are defined in terms of these diagrams. One reason we find the definition above to be more convenient is that boxes
should be thought of as ``colorless" for the purpose of many proofs.

Each relation above used at most two colors, just as the relations for $\sl_2$ (in the version of Definition \ref{defn:SBSBim-sl2-alt}) used at most one color. This follows the general
principle that relations amongst Soergel bimodules are governed by the relations which appear for finite parabolic subgroups.

Let us pause to note one consequence of these relations.

\begin{claim} For any two colors, we have an equality exemplified as below.
\begin{equation} \label{eq:C3} {
\labellist
\small\hair 2pt
 \pinlabel {$c_{r,b}$} [ ] at 94 37
\endlabellist
\centering
\ig{1}{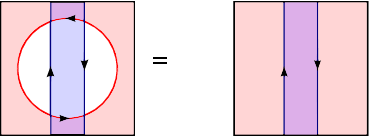}
} \end{equation} The scalar $c_{r,b}$ was given in Definition \ref{defn:roots3}.   \end{claim}

\begin{proof} We summarize the proof of \cite[Claim 6.5]{EDihedral}. By applying \eqref{singR2nonoriented1} to the right half of the LHS, and then sliding polynomials and using
\eqref{singR2}, one obtains the endomorphism of $rpr$ given in Sweedler notation by $\pa_r \pa_b \D^r_{p (1)} \ot \D^r_{p (2)} = c_{r,b} \pa^r_p \D^r_{p (1)} \ot \D^r_{p (2)} \in R^r
\ot_{R^p} R^r$. For any basis of $R^r$ over $R^p$, only a single term can survive the application of $\pa^r_p$ for degree reasons, and this is the term which pairs against the unit in the
dual basis. Thus $\pa^r_p \D^r_{p (1)} \ot \D^r_{p (2)} = 1 \ot 1$. \end{proof}

Now we define the singular Soergel $2$-category.

\begin{defn} \label{defn:ssbimtrue} Let $\SSBim_q = \SSBim_q(\sl_3)$ denote the following full sub-$2$-category of $\textbf{Cat}$. The objects are proper subsets of $S$, where $I
\subsetneq S$ is identified with the category of graded $R^I$-modules (using the notation $R^\mt = R$). The $1$-morphisms are summands of (grading shifts of) iterated induction and
restriction functors, between $R^I$-modules and $R^J$-modules for $I \subset J$. The $2$-morphisms are all natural transformations (i.e. homogeneous bimodule maps). \end{defn}

There is a $2$-functor $\GC$ from $\SBSBim_q$ to $\SSBim_q$ which is the identity on objects, sending downwards arrows to restriction functors (with grading shifts) and upwards arrows to
induction functors. The cups and caps are sent to units and counits of adjunction. An upwards crossing is sent to the natural isomorphism between iterated induction functors. The grading
shift on the restriction functor from $R^I$ to $R^J$ is $\ell(J)-\ell(I)$, where $\ell(p)=\ell(g)=\ell(o)=3$, $\ell(r)=\ell(b)=\ell(y)=1$, and $\ell(\mt)=0$; this is the length of the
longest element in the corresponding parabolic subgroup of the affine Weyl group.

\begin{claim} \label{claim:sl3ssbim} Suppose that $\Bbbk = \Q(q)$, or that $q=1$ and $\Bbbk$ is a field of characteristic $\ne 2,3$. Then $\GC$ is well-defined and fully faithful.
\end{claim}

The proof of this fact is a long algorithmic computation, and can be found in the appendix.

\subsection{The equivalence}
\label{subsec-equiv-sl3}

\begin{defn} Let $\mSBSBim_q$ denote the full sub-$2$-category of $\SBSBim_q$ defined as follows. The objects are $\{o,g,p\}$. The $1$-morphisms are generated by $IsJ$ whenever $I,J$ are
doubletons with $I \cap J = s$. The $2$-morphisms are full. In other words, our diagrams may contain any color, but the right-hand and left-hand colors must be in $\{o,g,p\}$, and the
bottom and top boundaries do not contain the color white. \end{defn}

\begin{defn} \label{def3fundtosing} Let $\FC$ be the $2$-functor $\FundOmq \to \mSBSBim_q$ defined as follows. On objects, it acts as the identity. On $1$-morphisms, it sends $go \mapsto g y o$, $og \mapsto o y g$ and so forth. On $2$-morphisms, it acts as below. \igc{1}{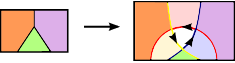} \end{defn}

\begin{claim} The $2$-functor $\FC$ is well-defined, and its image consists of degree $0$ maps. \end{claim}

\begin{proof} The isotopy relations in $\FundOmq$ follow from the isotopy relations in $\mSBSBim_q$. We check the remaining relations for a specific coloration, but the computation for
other colorations is almost identical.

Apply $\FC$ to the LHS of the circle relation \eqref{eq:circle3} to obtain a clockwise circle inside a anticlockwise circle. By \eqref{bigcircle} and \eqref{bigdemazure}, the result is
multiplication by $\pa^r_p(\mu^r_o)$. By \eqref{howtogetq3}, this is $[3]$.

Apply $\FC$ to the LHS of the bigon relation \eqref{eq:bigon}. \igc{1.5}{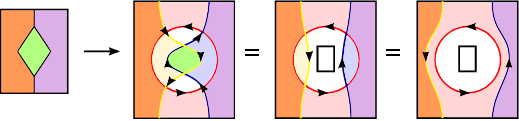} The first equality arises from \eqref{singR2nonoriented2}, and the second from \eqref{singR2}, while the polynomial in the box is $\mu_g^{b,y}$. Evaluating the red circle with \eqref{demazureis}, one obtains $\pa_r(\mu_g^{b,y})$, which
by \eqref{howtogetq2} is $-[2]$ as desired.

Apply $\FC$ to the LHS of the square relation \eqref{eq:square}.

{
\labellist
\small\hair 2pt
 \pinlabel {$c_{y,b}$} [ ] at 174 38
 \pinlabel {$c_{b,y}$} [ ] at 273 37
\endlabellist
\centering
\igc{1}{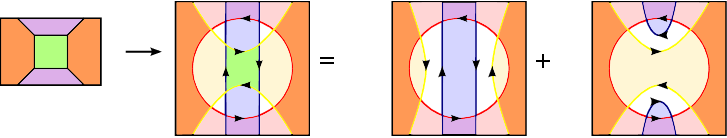}
}

The equality arises from relation \eqref{squaresbsb}. By applying \eqref{singR2} several times, this becomes

{
\labellist
\small\hair 2pt
 \pinlabel {$c_{y,b}$} [ ] at -10 36
 \pinlabel {$c_{b,y}$} [ ] at 85 37
\endlabellist
\centering
\igc{1}{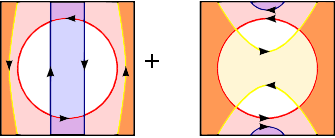}
}

after which we apply \eqref{eq:C3} to each side to obtain

{
\labellist
\small\hair 2pt
 \pinlabel {$c_{y,b} c_{r,b}$} [ ] at -20 36
 \pinlabel {$c_{b,y} c_{r,y}$} [ ] at 92 37
\endlabellist
\centering
\igc{1}{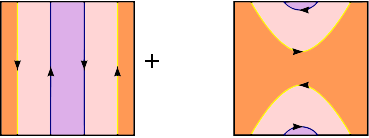}
}

Now the result follows from \eqref{howtoget1}. \end{proof}

\begin{remark} Most people would agree that representation theory is the ``easy" side of the geometric Satake equivalence, while perverse sheaves or Soergel bimodules are the ``hard" side,
at least in terms of the requisite background. Computationally, however, we like to argue that $\SBSBim_q$ is in some sense ``easier" to work with than $\Fund_q$, because the relations are
more local. For example, the generating trivalent vertex is sent by $\FC$ to a composition of several simpler maps, and the bigon and square relations each follow (after applying $\FC$)
from a sequence of relations on subdiagrams. \end{remark}

\begin{thm} \label{thm:sl3main} Suppose that $3$ is invertible in $\Bbbk$. Then the $2$-functor $\FC$ is a degree-zero equivalence. \end{thm}

\begin{proof} This proof follows a similar diagrammatic argument to the $\sl_2$ case, but is quite a bit more complicated. It can be found in the appendix. \end{proof}

\section{Representations and $\sl_n$-webs}
\label{sec-webs}

The goal of this chapter is to give some general introduction to the $2$-category of representations of a complex semisimple Lie algebra $\gg$. The focus will be on the algebraic,
generators and relations approach. We feel as though the literature is lacking an introduction of this form, though another similarly-focused introduction can be found in \cite{CKM}.

When we eventually apply this technology to geometric Satake, we will be applying it to the case of the Langlands dual $\gg^\vee$. However, in the interest of making this chapter more
readable, we stick with $\gg$ until \S\ref{subsec-om-versions}.

\subsection{Some algebraic philosophy}
\label{subsec-webs-philosophy}

Let $\Rep = \Rep(\gg)$ denote the monoidal category of finite dimensional representations of $\gg$, and let $\Rep_q = \Rep(U_q(\gg))$ denote its $q$-deformation.

The category $\Rep$ is semisimple, which implies that for any two irreducibles $V_\l$ and $V_\mu$ one has $\Hom(V_\l,V_\mu) = \delta_{\l \mu} \CM$. In other words, describing morphisms
between arbitrary objects (given a decomposition into irreducibles) is a trivial task. A semisimple category over $\CM$ is determined up to equivalence by its Grothendieck group,
reinforcing the fact that the morphisms do not play an essential role.

However, $\Rep$ has a monoidal structure, and a semisimple monoidal category is not determined by its Grothendieck ring (see, for instance, \cite[Problem 1.42.8 and
following]{EtingofTensorNotes}). As such, we should study its morphisms in more detail.

\begin{goal} \label{goal:total} Describe the morphism spaces $\Hom(V_{\l_1} \ot \cdots \ot V_{\l_d}, V_{\mu_1} \ot \cdots \ot V_{\mu_{e}})$ for any sequences of irreducible
representations, along with the structures of vertical (i.e. usual) and horizontal (i.e. monoidal) composition. \end{goal}

This appears to be an incredibly difficult problem, and it remains open for almost all $\gg$. A more tractable approach is to study $\Fund = \Fund_\gg$, the full monoidal subcategory whose
objects are tensor products $V_{\o_{i_1}} \ot \cdots \ot V_{\o_{i_d}}$ of fundamental representations. Every irreducible is a summand of such a tensor product, so that the \emph{idempotent
closure} or \emph{Karoubi envelope} $\Kar(\Fund)$ is isomorphic to $\Rep$. An introduction to Karoubi envelopes can be found in \cite{BarMor}, or on Wikipedia. Here are three main goals in
the algebraic study of $\Fund$ and $\Rep$, with a discussion of the known results.

\begin{goal} \label{goal:Fund} Describe the monoidal category $\Fund$ by generators and relations. \end{goal}

The language of description tends to be planar diagrammatics. Of course, the presentation will depend on a choice of generating morphisms (just as our category $\Fund$ came from a choice of
generating objects), but choices can be made so that the relations take on a nice form.

\begin{itemize}
\item When $\gg = \sl_2$, $\Fund$ is isomorphic to the Temperley-Lieb category \cite{TemLie}, which was
given a diagrammatic interpretation by Kauffman \cite{Kau}. The diagrams involved are called crossingless matchings.
\item When $\gg$ has rank 2, a presentation for $\Fund$ was found by Kuperberg \cite{Kup}. His term for the diagrams, \emph{webs}, has become standard, so that crossingless matchings can also be called $\sl_2$-webs. 
\item When $\gg = \sl_n$, a conjectural presentation for $\Fund$ due to Morrison \cite{morrison2007diagrammatic} (and Kim \cite{kim2003graphical} for $n=4$) was recently proven by Cautis-Kamnitzer-Morrison \cite{CKM}, in terms of $\sl_n$-\emph{webs}.
\item Outside of these cases, a presentation for $\Fund$ is unknown.
\end{itemize}

In this paper we rely only on the description of $\Fund(\sl_n)$ by $\sl_n$-webs, as well as an analogous description of a nice subcategory of Soergel bimodules. However, to motivate the
introduction of $\Fund$, we should at least explain how one would complete this approach into a solution to the original problem, Goal \ref{goal:total}.

Let $\l = \sum a_i \o_i$ be a decomposition of a dominant integral weight $\l$ as a sum of fundamental weights $\o_i$, so that $a_i \in \N$. Let $V = V_{\o_{i_1}} \ot \cdots \ot
V_{\o_{i_d}}$ be some tensor product where each fundamental representation $V_{\o_i}$ appears $a_i$ times. Then $V_\l$ is a summand of $V$ with multiplicity one, so that there is some
canonical idempotent $e_V \in \End(V)$ which projects to $V_\l$. In the literature, such idempotents are often called \emph{clasps}.

\begin{goal} \label{goal:clasps} Using the algebraic description of $\Fund$ from the previous goal, give a (possibly recursive) formula for $e_V$. \end{goal}

\begin{itemize}
\item When $\gg = \sl_2$, the clasps are known as Jones-Wenzl projectors \cite{Jon1,Wenzl}, and several recursive formulas \cite{Wenzl, BFK} and a closed formula \cite{Mor} are known.
\item When $\gg$ has rank 2, a recursive formula for clasps was found by Kim \cite{kim2003graphical}.
\item The author \cite{ELLCC} has recently conjectured a recursive formula for clasps when $\gg = \sl_n$, proven to hold for $n \le 4$.
\item Outside of these cases, nothing is known.
\end{itemize}

Goal \ref{goal:clasps} is too naive, so we replace it by Goal \ref{goal:intertwiners} below (though it seems unlikely that one could solve the former without the latter). Consider the
finite set $\XM_\l$ consisting of all tensor products of fundamental representations whose fundamental weights add up to $\l$ (they are all isomorphic). For any two $V, W \in \XM_\l$,
there is a $1$-dimensional space of maps $V \to W$ which factor through the common summand $V_\l$.

\begin{goal} \label{goal:intertwiners} Find a family $\phi_\l$, consisting of maps $\phi_{V,W} \co V \to W$ for each $V, W \in \XM_\l$, such that $\phi_{V,W}$ factors through $V_\l$, and
such that $\phi_{W,X} \circ \phi_{V,W} = \phi_{V,X}$ for all $V, W, X \in \XM_\l$. \end{goal}

Note that $\phi_{V,V}$ is none other than the clasp $e_V$. Together, the family $\phi_\l$ gives a canonical identification of $V_\l$ as a common summand of each object in the family
$\XM_\l$.

Given a presentation of a monoidal category and an intertwining family $\phi$ (such as $\phi_\l$), it is a trivial operation to produce a presentation of the partial Karoubi envelope,
where the image of $\phi$ is formally added as a new object. In the literature (e.g. \cite{KLMS}), such a planar calculus which describes morphisms in a partial Karoubi envelope is called
a \emph{thick calculus}, or perhaps a \emph{thicker calculus} (with the thickest calculus being the entire Karoubi envelope). Adding $\phi_\l$ for each dominant weight $\l$, we obtain a
solution to Goal \ref{goal:total}.

This need not be the end of the story, as the thick calculus may itself have a number of interesting formulas and consequences, and this may simplify the presentation. For example, the
$3j$ and $6j$ symbols and theta networks in $\sl_2$ representation theory are computations within the thick calculus. For more about these topics, see \cite{FSS} and the references it
contains.

One additional feature of the known diagrammatic presentations is that they are easily adapted to the study of quantum group representations. A $q$-deformation of the $\sl_n$-web algebras
discussed above will yield a description of $\Fund_q$, tensor products of fundamental representations of $U_q(\sl_n)$. When $q$ is specialized to a root of unity (with $U_q(\sl_n)$
denoting the Lusztig form), representations of $U_q(\gg)$ are no longer semisimple; $\Kar(\Fund_q)$ will no longer be equivalent to $\Rep_q$, being equivalent to the category of tilting
modules instead. By understanding which denominators in the formulas for clasps will vanish at a root of unity, one can address the behavior of these tilting modules in an algebraic
manner. Of course, this is already known using abstract representation theory, but knowing the clasps explicitly will still yield measurable gains. In the broader context of Soergel
bimodules, computing the denominators in analogous idempotents is of great interest in modular representation theory.

\subsection{Monoidal-graded categories}
\label{subsec-om-monoidal}

\begin{defn} (c.f. \cite{ENO}) Let $H$ be a finite abelian group. A monoidal category $\CC$ is $H$-\emph{graded-monoidal} if there are subcategories $\CC_{\xi}$ for $\xi \in H$, for which
$\CC = \oplus_{\xi \in H} \CC_{\xi}$, and for which $\CC_\xi \ot \CC_\nu \subset \CC_{\xi \nu}$. \end{defn}

\begin{defn} The \emph{fundamental group} $\pi(\gg)$ is defined as $\L_{\wt} / \L_{\rt}$, the quotient of the weight lattice by the root lattice. We also denote it $\Om$. \end{defn}

The category $\Rep = \Rep(\gg)$ is $\Om$-graded-monoidal (as is $\Rep_q$). For $\xi \in \Om$, let $\Rep_\xi$ denote those representations whose irreducible constituents all have highest
weights in the coset $\xi$. Clearly the tensor product obeys the group law. From this monoidal grading one can reconstruct the category of $G$ representations for any corresponding
algebraic group.

\begin{example} Representations of $\sl_2$ or $SL_2$ can be split into even and odd representations. These can be distinguished by the action of the subgroup $\{ \pm 1 \} \subset SL_2$ on
these representations. The subcategory of even representations is equivalent to the category of $\PM SL_2$-representations. \end{example}

\begin{remark} The adverb $H$-\textbf{graded} modifies the adjective \textbf{monoidal}. One should not confuse this with other uses of the term \textbf{graded}. For instance, a category is
\emph{graded} or \emph{$\Z$-graded} if (by conflating maps of all ``degrees") its Hom spaces can be enriched in graded vector spaces. Therefore, every (homogeneous) \textbf{morphism} has
an associated degree in $\Z$. On the other hand, every (indecomposable) \textbf{object} in an $H$-graded-monoidal category has an associated character in $H$. \end{remark}

\begin{constr} It is a trivial operation to take an $H$-graded-monoidal category and replace it with a $2$-category $\CG$ which encodes (almost) the same data. Let $S$ be a right
$H$-torsor. The objects of $\CG$ will be identified with $s \in S$, and for $\xi \in H$ the morphism category from $s$ to $s \xi$ is $\CC_{\xi}$. A $2$-category where the objects form
an $H$-torsor, equipped with a corresponding $H$ action which identifies $\Hom(s,t) = \Hom(s \xi,t \xi)$, we will call an \emph{$H$-$2$-category}, for lack of better terminology.
\end{constr}

\begin{remark} The data of an $H$-$2$-category $\CG$ together with a choice of $s \in S$ (corresponding to $1 \in H$) and the data of an $H$-graded-monoidal category $\CC$ are equivalent.
Philosophically speaking, the difference between them is analogous to the difference between an algebra and its regular representation. Because $\CC$ is monoidal, it acts on itself by
tensor product. Under the splitting $\CC = \oplus_{\xi \in H} \CC_{\xi}$, the action of $M \ot (\cdot)$ for any object $M$ can be expressed as an $H \times H$ matrix. Each entry of this
matrix is a corresponding 1-morphism in $\CG$. \end{remark}

Thus from $\Rep$ (resp. $\Rep_q$) we can construct a $\Om$-$2$-category, which we denote $\RepOm$ (resp. $\RepOmq$). It is simple to pass from a diagrammatic description of $\Rep(\gg)$
(say, by $\sl_n$-webs) to a diagrammatic description of $\RepOm$. One simply colors each planar region by an element of $\Om$, and imposes the rule that tensoring with $V_{\o_i}$ will
shift the color by $\o_i$. This is illustrated by the $\sl_2$ and $\sl_3$ cases in previous chapters.

The reason we switch from graded monoidal categories to $2$-categories is that it makes it easier to compare $\RepOmq$ with $\SSBim_q$, which is more naturally a $2$-category. When $q=1$,
the subcategory $\mSBSBim_q$ which is the target of the algebraic Satake equivalence is equivalent to a graded $\Om$-$2$-category, and thus arises from some $\Z$-graded
$\Om$-graded-monoidal category, though this is perhaps an unnatural way of viewing it.

\begin{remark} The monoidal category $\Rep$ has a symmetric structure. There is no notion of a symmetric structure on a $2$-category, but the corresponding structure on an
$\Om$-$2$-category is not hard to formulate. \end{remark}

\subsection{The fundamental group}
\label{subsec-om-versions}

We return to the setting of geometric Satake, which cares about representations of $\gg^\vee$. We provide several important facts about the fundamental group $\Om = \pi(\gg^\vee)$. It is a
finite abelian group, so that its Pontrjagin dual or character group $\Om^*$ is also a finite abelian group, non-canonically isomorphic to $\Om$.

Let $\G$ denote the Dynkin diagram of $\gg$, and $\G^\vee$ the Langlands dual Dynkin diagram. Write $G_\adj$ or $G_\sc$ for the adjoint or simply-connected algebraic groups for
$\G$. Let $\KC = \C((t))$ and $\OC = \C[[t]]$.

\begin{claim} The following finite abelian groups are all canonically isomorphic, and will all be denoted $\Om$.
\begin{itemize}
	\item $\pi(\gg^\vee)$.
	\item $\pi(\gg)^*$.
	\item $Z(G_{\sc})$, the center of $G_{\sc}$.
	\item $\pi_1(G_{\adj})$.
	\item $\pi_0(G_{\adj}(\KC))$, the component group of $G_{\adj}(\KC)$.
\end{itemize}
\end{claim}

Let $\tG$ denote the affine Dynkin diagram of $\G$, and let $0$ denote the affine vertex (i.e. choose an embedding $\G \into \tG$). Call a vertex $v \in \tG$
\emph{removable} if $\tG \setminus v$ is isomorphic to $\G$, and let $\trmv$ denote the set of removable vertices. They form a single orbit under the automorphism group of
$\tG$. Let $\rmv = \trmv \setminus 0 \subset \G$. Identifying the vertices of $\G^\vee$ and $\G$, we can view $\rmv$ as a subset of $\G^\vee$. It is known that the fundamental
representations of $\gg^\vee$ associated to vertices in $\rmv$ are precisely the \emph{miniscule} fundamental representations, meaning that their weights form a single Weyl group orbit.

\begin{claim} The set $\trmv$ is in canonical bijection with $\Om$. That is, it is an $\Om$-torsor with a distinguished element $0$. \end{claim}

Essentially, the statement is that the miniscule fundamental weights, together with the zero weight, enumerate a list of representatives in $\L^\vee_{\wt}/\L^\vee_{\rt} = \Om$.

\begin{remark} Here, Langlands duality is essential. The miniscule fundamental weight in type $B$ is precisely the removable vertex from affine type $C$, and vice versa. \end{remark}

These facts are quite standard, see for instance \cite[Chapter 2]{wu2008miniscule}.

\section{Singular Soergel Bimodules}
\label{sec-ssbim}

The goal of this chapter is to provide background material on algebraic Soergel theory. For the connections with geometry, see \S\ref{sec-reformulate}. This background is not essential to
the main theorems.

The exotic realizations of affine $\sl_n$ for $n > 2$ are examples of ``odd-unbalanced" realizations, as defined in \cite{EDihedral}, while the ordinary realization of affine $\sl_n$ is
odd-balanced. Being odd-unbalanced adds an additional layer of technicality which readers uninterested in the $q$-deformation of geometric Satake can ignore. Because of this, it may be
better for readers to get their first introduction to singular Soergel bimodules elsewhere, such as in \cite{WilSSB}. This paper is the first to use odd-unbalanced realizations in an
essential way, which is why we devote a great deal of effort to explaining these technicalities.

\subsection{Overview}
\label{subsec-summary}

Here is a big picture to keep in mind while reading the nitty-gritty. 

Given a (reflection faithful) representation $\hg$ of a Coxeter group $W$, one can consider the polynomial ring $R = \OC(\hg)$. This is a graded ring whose linear terms $\hg^*$ will be
placed in degree $2$, equipped with a natural action of $W$. The letters $I,J$ will designate subsets of the simple reflections $S$ of $W$. For any parabolic subgroup $W_I \subset W$ there
is a subring of invariants $R^I \subset R$. Under mild assumptions, the ring extension $R^s \subset R$ is graded Frobenius, meaning that induction and restriction functors are biadjoint
(up to a shift).

\begin{defn} Let $\BSBim = \BSBim_q$, the monoidal category of \emph{Bott-Samelson bimodules}, be the full subcategory of $(R,R)$-bimodules whose $1$-morphisms are generated by \[B_s = R
\ot_{R^s} R (1)\] for each simple reflection $s \in S$, the composition of restriction with induction and a shift. The grading shift places $1 \ot 1$ in degree $-1$. \end{defn}

\begin{defn} The monoidal category $\SBim$ of \emph{Soergel bimodules} is the Karoubi envelope of $\BSBim$. \end{defn}

Thus to define $\SBim$ one only needs the data of a reflection representation $\hg$ satisfying some conditions. Under some additional conditions, Soergel \cite{Soe5} has proven that
$\SBim$ has a number of desirable properties: it categorifies the Hecke algebra, and the Hecke algebra is equipped with a pairing which encodes the graded dimension of morphism spaces in
$\SBim$.

However, our goal is not just to define Soergel bimodules but to find an algebraic presentation for them. One must first choose a set of generating morphisms. As it turns out, this arises
from a choice of Frobenius structure on the extensions $R^s \subset R$ for $s \in S$, which boils down to an explicit choice of simple roots in $\hg^*$. A reflection representation paired
with a choice of simple roots is (essentially) what we call a \emph{simple realization} below, the word \emph{simple} indicating that is corresponds to a choice of simple roots. A simple
realization can roughly be encoded in a Cartan matrix. From a simple realization one can define a diagrammatic category by generators and relations which is equivalent to $\BSBim$.
Moreover, for a slightly more general notion of a simple realization (allowing for non-faithful representations and more) one can still define this diagrammatic category and prove the
categorification results of Soergel, even when it is no longer equivalent to $\BSBim$. This is shown in \cite{EWGR4SB}.

Now return to the original setting of a (reflection faithful) representation $\hg$. Under certain assumptions, whenever $W_I \subset W_J$ are two finite parabolic subgroups, the ring
extension $R^J \subset R^I$ will be graded Frobenius (this is an upgraded version of Chevalley's Theorem \cite[Chapter 3]{Humphreys}). Let $\ell(I)$ denote the length of the longest
element of $W_I$.

\begin{defn} Let $\SBSBim = \SBSBim_{\textrm{alg}}$, the $2$-category of \emph{singular Bott-Samelson bimodules}, be defined as follows. The objects are the rings $R^I$ for finite
parabolic subgroups $W_I$. The morphism category $\Hom(R^I,R^J)$ will be a full subcategory of $(R^J,R^I)$-bimodules. The $1$-morphisms in $\SBSBim$ are generated by the induction bimodule
${}_{R^I} R^I _{R^J}$ and the restriction bimodule ${}_{R^J} R^I _{R^I} (\ell(J)-\ell(I))$ for the extensions $R^J \subset R^I$ whenever $I \subset J$. \end{defn}

\begin{defn} The $2$-category $\SSBim$ of \emph{singular Soergel bimodules} is the Karoubi envelope of $\SBSBim$. \end{defn}

Thus to define $\SSBim$ one only needs the data of a representation $\hg$, satisfying some conditions. Under additional conditions, Williamson \cite{WilSSB} has proven that $\SSBim$ has a
number of desirable properties: it categorifies the Hecke algebroid, and the Hecke algebroid is equipped with a pairing which determines the graded rank of 2-morphism spaces.

To find an algebraic presentation for $\SBSBim$ one must choose a set of generating 2-morphisms. Again, this arises from a choice of Frobenius structures on every extension $R^J \subset
R^I$ together with some consistency condition on these structures. We call this data a \emph{Frobenius realization} (Definition \ref{defn:frobmaps}). One way to specify a Frobenius
realization is to make an explicit choice of positive roots in $\hg^*$, which we call a \emph{root realization} (Definition \ref{defn:roots3}). From a Frobenius realization one can
define a diagrammatic category by generators and relations which is equivalent to $\SBSBim$. As above, one can generalize the notion of a Frobenius or root realization, and one should be
able to define $2$-categories which possess the categorification properties proven by Williamson, even without being equivalent to $\SBSBim$. This will eventually be proven in type $A$ in
\cite{EWGR4SSB}. It was proven for dihedral groups in \cite{EDihedral}, and will be proven for affine $\sl_3$ in the appendix.

Let $\a_s$ denote the simple root attached to $s \in S$. In familiar notions of root systems (as in, say, Humphreys' book \cite{Humphreys}), a choice of simple roots determines a choice of
positive roots. One definition is that the positive roots are the elements of the set $\{ s_1 s_2 \ldots s_{d-1} (\a_{s_d}) \}$ for each reduced expression $\un{w} = s_1 \ldots s_d$. This
set has redundancies: for example, when $m_{st} = 3$, one has $s(\a_t) = t(\a_s)$. However, for more general realizations, the previously redundant descriptions of the same positive root
may differ by a scalar! For instance, using the notation from \S\ref{subsec-polynomials-sl3}, one has $m_{by}=3$ but $y(\a_b) = q b(\a_y)$. This kind of ambiguity can occur whenever the
Cartan matrix of the realization is odd-unbalanced, in which case the choice of positive roots is a piece of additional data. In a precisely analogous computation, the Demazure operators
attached to simple roots will be used to define the trace maps for all Frobenius extensions $R^J \subset R^I$, but they only satisfy the braid relations when the realization is
odd-balanced.

There is not a great deal of complication involved - just keeping track of some scalars which measure the failure of roots to be well-defined! The notion of a root realization is developed
as a natural way to keep track of these scalars. Remember also that the $2$-categories themselves only depend on the representation $\hg$, while the presentation depended on the root
realization.

The sections to follow will explore these concepts in more detail, with constant reference to the exotic realization of affine $\sl_n$. The one remaining subtlety is that, unlike the
familiar affine Cartan matrix, the exotic affine Cartan matrix is non-degenerate. This may lead to possible confusion when comparing the exotic and familiar realizations, which we do our
best to forestall. This motivates the discussion of dual realizations below.

\subsection{Simple realizations}
\label{subsec-simplerealize}

Fix a Coxeter group $W$ with a set of simple reflections $S$, and let $m_{st}$ be the order of $st$, for $s,t \in S$. Let $\Bbbk$ be a commutative ring. We simplify the discussion somewhat
by assuming that $\Bbbk$ has no $2$-torsion. Typical choices for $\Bbbk$ are $\C$, $\Q(q)$, and $\Z[\z_n]$ for a primitive $n$-th root of unity $\z_n$.

\begin{defn} \label{defn:realization} A \emph{(simple) realization} of $(W,S)$ over $\Bbbk$ is a free, finite rank $\Bbbk$-module $\hg$, together with subsets $\D^\vee = \{ \a_s^\vee \; | \; s \in S\} \subset \hg$ and $\D = \{ \a_s \; | \; s \in S\} \subset \hg^* = \Hom_{\Bbbk}(\hg,\Bbbk)$, called \emph{simple coroots} and \emph{simple roots} respectively, satisfying:
\begin{enumerate}
\item $\langle \alpha_s^\vee, \alpha_s \rangle = 2$ for all $s \in S$;
\item the assignment \[s(v) \define v - \langle v, \alpha_s\rangle \alpha_s^{\vee}\] for all $v \in \hg$ yields a representation of $W$.
\end{enumerate}
We will often refer to $\hg$ as a realization, however the choice of simple roots and coroots is always implicit. \end{defn}

Note that $\hg^*$ is equipped with a contragredient action of $W$, for which \[s(f) = f - \langle \a_s^\vee, f \rangle \a_s. \]

\begin{defn} Given a realization $(\hg,\hg^*,\D^\vee,\D)$, the \emph{dual realization} is $(\hg^*,\hg,\D,\D^\vee)$. This makes sense, since $\hg^{**} \cong \hg$ canonically. \end{defn}

\begin{defn} \label{defn:cartan} The \emph{Cartan matrix} associated to a realization $\hg$ is the $S \times S$ matrix $A$ valued in $\Bbbk$, with entries $a_{s,t} = \langle a_s^\vee, \a_t
\rangle$. \end{defn}

Note that we did not assume that either $\D$ or $\D^\vee$ is a basis, or is even linearly independent. When $\D^\vee$ is a basis for $\hg$, the realization is determined by its Cartan
matrix. However, ``enlarging" the representation beyond the span of $\D$ may affect the image of $\a_s$ or $\a_s^\vee$ in $\Bbbk$.

\begin{example} Let $W$ be a Weyl group equipped with a root system. Then the familiar reflection representation (with symmetric Cartan matrix, as in \cite{Humphreys}) is a self-dual
realization of $W$ over $\R$. \end{example}

\begin{example} \label{ex:stdaffine} Let $W$ be an affine Weyl group, with $n+1$ simple reflections. It has a familiar realization over $\RM$ of rank $n+1$ where $\D$ forms a basis, while
$\D^\vee$ is linearly dependent. It has a dual realization of rank $n+1$ where $\D^\vee$ forms a basis, and $\D$ is linearly dependent. It has a self-dual realization of rank $n+2$ where
both $\D$ and $\D^\vee$ are linearly independent, but neither is a basis. \end{example}

\begin{claim} Suppose that the Cartan matrix has invertible (not just non-zero) determinant in $\Bbbk$. Then $\D$ and $\D^\vee$ are each linearly independent. As $W$ representations one
has $\hg \cong \Bbbk \cdot \D^\vee \oplus \TM$ for some trivial representation $\TM$ in the kernel of $\D$. Similarly, one has $\hg^* \cong \Bbbk \cdot \D \oplus \TM^*$. \end{claim}

\begin{proof} The fact that $\D$ and $\D^\vee$ are linearly independent is obvious. Let $v$ be a vector in $\hg$ which is not in $\Bbbk \cdot \D^\vee$. Using linear algebra, some
vector $x \in \Bbbk \cdot \D^\vee$ satisfies $\a_s(x) = \a_s(v)$ for all $s \in S$. Then $v-x$ is in the kernel of $\D$, so that $W$ acts trivially on $v-x$. In this fashion, one can
replace any basis for $\hg$ extending $\D^\vee$ by a basis where the remaining elements are in the kernel of $\D$. \end{proof}

It is obvious how one equips the base change $\hg \ot_{\Bbbk} \Bbbk'$ with the structure of a realization over $\Bbbk'$. If $\D^\vee$ spans a free module of rank $|S|$ in $\hg$, then it
does the same in $\hg \ot_{\Bbbk} \Bbbk'$. However, the same can not be said for $\D$. For any Cartan matrix over a domain $\Bbbk$ with non-zero non-invertible determinant $D$, base change
to $\Bbbk/D$ will lower the rank of the span of $\D$. Therefore, base change does not commute with duality.

\subsection{The exotic realization}
\label{subsec-exoticrealization}

In this section we let $W$ be the affine Weyl group in type $\tilde{A_n}$, and let $\Bbbk = \Z[q,q^{-1}]$. The finite Weyl group $W_\fin$ will be embedded inside $W$, generated by the
first $n$ simple reflections.

\begin{defn} We define the \emph{exotic affine Cartan matrix of $\sl_{n+1}$}. For $n\ge 2$ it has a form which is exemplified by the case $n=4$:

\[A = \left( \begin{array}{ccccc} 2 & -1 & 0 & 0 & -q^{-1} \\ -1 & 2 & -1 & 0 & 0 \\ 0 & -1 & 2 & -1 & 0 \\ 0 & 0 & -1 & 2 & -q \\ -q
& 0 & 0 & -q^{-1} & 2 \end{array} \right).\]

Exotic affine $\sl_2$ on the other hand is given by a symmetric Cartan matrix:

\[A = \left( \begin{array}{cc} 2 & -(q+q^{-1}) \\ -(q+q^{-1}) & 2 \end{array} \right). \]

We call any realization over $\Zqq$ with this matrix an \emph{exotic realization} of $W$. \end{defn}

Obviously, setting $q=1$ yields the usual affine Cartan matrix for $\tilde{A_n}$.

\begin{exercise} Whenever $a_{st} a_{ts}=1$, one has $sts=tst$. Therefore, this matrix gives a realization of $W$. \end{exercise}

Thus whenever $n \ge 2$, specialization of $q$ will not change the order of $s_i s_j$ for any simple reflections, so that the realization is dihedrally faithful (as in
\cite{EDihedral}).

\begin{exercise} The exotic Cartan matrix has determinant $2 - q^2 - q^{-2}$, which is non-zero and non-invertible. Therefore, in any base change of the exotic realization, $\D^\vee$ will
be linearly independent. \end{exercise}

\begin{exercise} For the exotic realization of rank $n+1$ where $\D^\vee$ is a basis, base change to $q^2 = 1$ will make $\D$ linearly dependent. This is the realization for which the
action of $W$ on $\hg^*$, after specialization to $q=1$, is the reflection representation of $W$ defined in \cite{Humphreys}. In particular, an exotic realization is faithful. There are
also exotic realizations for which base change to $q=1$ will yield the other realizations mentioned in Example \ref{ex:stdaffine}. \end{exercise}

Number the simple roots so that $\D = \D_{\fin} \cup \{\a_0\}$, where $\D_{\fin}= \{\a_1,\ldots, \a_n\}$ has the usual finite Cartan matrix and $\a_0$ is the affine simple root. Let $s_i$
denote the reflection corresponding to $\a_i$. Let $\b = \sum_{i=1}^n \a_i$ and $\b^\vee = \sum_{i=1}^n \a_i^\vee$, the highest finite root and coroot. Let $t$ denote reflection across
$\b$, which can be written as $t = s_1 s_2 \cdots s_{n-1} s_n s_{n-1} \cdots s_2 s_1$. Let $[n]$ denote the $n$-th quantum number in $\Zqq$.

\begin{exercise} The subgroup $W'$ generated by $s_0$ and $t$ inside $W$ is the infinite dihedral group. Let $\D' = \{\beta,\a_0\}$ and $(\D^\vee)' = \{\b^\vee,\a_0^\vee\}$. Then
$(\hg,\hg^*,\D',(\D^\vee)')$ gives a realization of $W'$, having as its Cartan matrix the exotic affine Cartan matrix of $\sl_2$. \end{exercise}

In \cite{EDihedral} the exotic realization of $\sl_2$ was studied in detail. It was shown that the realization is faithful unless $q^2$ is a root of unity. If $q^2$ is a primitive $m$-th
root of unity then the kernel of the action is generated by $(uv)^m$ (for simple reflections $u,v$).

\begin{prop} After base change, an exotic realization of $\sl_{n+1}$ is still faithful unless $q^2$ is a primitive $m$-th root of unity, in which case the kernel $K$ is generated by $(s_0
t)^m$. \end{prop}

Presumably there is an elegant proof of this fact, but we provide a brute force proof.

\begin{proof} We assume that $q^2$ is a primitive $m$-th root of unity in the rest of the proof; it is an exercise to modify this into a proof of faithfulness when $q^2$ is not a root of
unity. Our first step is to examine $H = K \cap \L_{\rt}$, those translations within the kernel. Clearly $H$ is closed under the action of $W_\fin$. We wish to show that $H = m\L_{\rt}$.  

The element $(s_0 t)^m$ is translation by the $m$-th multiple of the highest root. Since this translation is in $H$, the entirety of $m \L_{\rt}$ is in $H$. Moreover, $(s_0 t)^k \in H$ if
and only if $k$ is a multiple of $m$ (by the dihedral result above), meaning that a $k$-th multiple of a root is in $H$ if and only if $k$ is a multiple of $m$.

If $\l \in H$ then $s_i(\l) - \l = \langle \l, \a_i \rangle \a_i \in H$, which is a multiple of a root. Therefore $\langle \l, \a_i \rangle$ is a multiple of $m$ for all $i$, meaning that
$\l \in m \L_{\wt}$. Thus $m\L_{\rt} \subset H \subset m\L_{\wt}$. Since $m\L_{\wt} / m\L_{\rt} \cong \Z/(n+1)\Z$, there are not many choices for $H$: one has $H = H_k$ for some $k$
dividing $n+1$, where $H_k$ is generated by $m\L_{\rt}$ and $km\L_\wt$. Moreover, $H \subset \L_\rt$, which places further restrictions. We seek to show that $k=n+1$, in which case
$km\L_\wt$ is contained inside $m\L_\rt$.

Now consider the element $s_1 s_2 \ldots s_0 = x$. A simple computation shows that $x$ preserves the subspace spanned by $\{ \a_2, \ldots \a_n\} \cup \{q^{-1} \a_1 + \a_0 = s_1(\a_n) \}$,
and that it acts by a rescaled permutation matrix. From this it is easy to compute that the order of $x$ is precisely $mn$. Note that $x^n$ is translation by an $(n+1)$-st multiple of a
fundamental weight which, because $\Om = \Z/(n+1) \Z$, is in $\L_{\rt}$. Let $\o$ be the fundamental weight with $(n+1)\o = x^n$.

Let $d$ be the greatest common divisor of $n+1$ and $m$, so that $n+1 = da$ and $m=db$. The fact that $km\L_\wt \subset \L_\rt$ implies that $n$ divides $km$, meaning that $k=al$ for some
$l$ dividing $d$. However, $K$ contains $km(\o) = bl(n+1)\o = (x^n)^{bl}$. By previous computation, this is impossible unless $l=d$ and $k=n+1$. Thus we have finally proven that $H =
m\L_\rt$.

Our next step is to examine $J$, the image of $K$ under the quotient map $W \to W_\fin$ which kills all translations. Clearly $J$ is a normal subgroup of $S_{n+1}$, which means that it is
trivial or it contains the alternating subgroup $A_{n+1}$ (or $n=1$). If we show that $J$ is trivial then we have proven the proposition. Thus we aim for a contradiction when $s_1 s_2 \in
J$ (or $s_1 \in J$ for $n=1$). We sketch the rest of the proof.

Let $\l$ denote translation by $\l \in \L_\rt$, and suppose that $w \l \in K$ for $w \in W_\fin$. By applying $H$ we may freely shift $\l$ by an element of $m \L_\rt$. Clearly $\l \notin m
\L_\rt$, because $W_\fin$ acts faithfully. By a fairly straightforward argument using convolution and multiplication, one can show that $w(\l)-\l \in K$, which implies by the result above
that $w(\l)-\l \in m\L_\rt$. If this is the case for $w = s_1$ or $w = s_1 s_2$, it is a similar argument to show that $s_i(\l)-\l \in m\L_\rt$ for all $i \le n$, or in other words, that
$\l \in m\L_\wt$. Letting $d$ be the greatest common divisor of $n+1$ and $m$, and $a,b$ as above, this implies that $\l \in b \L_\rt \cap m \L_\wt$. After some manipulation we may assume
that $\l = (x^n)^{bl}$ for some $l$ dividing $d$.

Now the result follows from computing $s_1 x^{nk}$ and $s_1 s_2 x^{nk}$ on the subbasis $\{ \a_2, \ldots \a_n\} \cup \{q^{-1} \a_1 + \a_0 = s_1(\a_n) \}$, and observing that these
operators are never trivial. \end{proof}

We conclude this section with a discussion of where the exotic representations of $W$ can be found in the literature, and some history. This was explained to me by Lusztig.

There is another ``deformation of the Cartan matrix" which may be more familiar to readers: one with the usual (integral) non-diagonal entries, but replacing the diagonal entries $2$ with
$v+v^{-1}$. (This is a different formal parameter $v$, not the same as $q$!) Warning: this is not a Cartan matrix, in the sense we have defined above. This deformation, originally found by
Killing, was used by Kilmoyer in his thesis \cite{kilmoyerThesis} (see also \cite{kilmoyerEtc}) in finite type to define a representation of the Hecke algebra (with parameter $v$) which
deforms the usual reflection representation. Analogously, the $v$-deformed $q$-exotic Cartan matrix provides a one-parameter family of affine Hecke algebra representations.

Kazhdan and Lusztig provide a definition of $W$-graphs \cite{KaLu1}, giving $v$-deformations of other $W$ representations. In \cite{LusSqint}, Lusztig gives
a number of $W$-graphs for affine Weyl groups, using a special left cell; in affine types $A$ and $C$, these $W$-graphs can be observed to have periodicity properties. In
\cite{LusPeriodic}, it is explicitly described how to take a periodic $W$-graph and obtain a one-parameter family of Hecke representations. In type $A$, the corresponding family of affine
Hecke representations, after passing to $v=1$, is precisely the exotic family of representations of $W$.

However, our parametrization of this family by the variable $q$ appears to be new. Accordingly, study of what happens when $q$ is specialized to a root of unity also appears to be new.

\begin{remark} There is also a one-parameter family of reflection representations in affine type $C$. One might pray that a $q$-parametrization might exist which gives a quantum Soergel
Satake equivalence. \end{remark}

%
%
%
%

\subsection{Frobenius extensions}
\label{subsec-frobenius}

\begin{defn} A \emph{graded (commutative) Frobenius extension of degree $d$} is an extension of graded commutative rings $A \subset B$ equipped with a \emph{trace map} $\pa \co B \to A$
which is $A$-bilinear and has degree $-d$. There must exist bases $\{b_i\}$ and $\{b_j^*\}$ of $B$ as a finite rank free $A$-module, satisfying $\pa(b_i b_j^*) = \d_{i,j}$. \end{defn}

The choice of dual bases is not part of the structure of a Frobenius extension, but the choice of trace map is. For any ring extension there are two canonical bimodule maps: inclusion $i
\co A \to B(d)$ of $A$-bimodules, having degree $-d$, and multiplication $m \co B \ot_A B(d) \to B$ of $B$-bimodules, having degree $+d$. For a Frobenius extension there are two additional
bimodule maps: trace $\pa \co B(d) \to A$ of $A$-bimodules, having degree $-d$, and comultiplication $\D \co B \to B \ot_A B(d)$ of $B$-bimodules, having degree $+d$. The comultiplication
map satisfies $\D(1) = \D^B_A = \sum b_i \ot b_i^*$, this element being independent of the choice of dual bases.

Note that for any invertible scalar $\l$, $\l \pa$ will also define a trace map, producing a different Frobenius extension structure. One way to ``measure" the Frobenius structure is by the
invariant $\mu^B_A = m(\D(1)) \in B$, the \emph{product-coproduct element}, which has degree $2d$. Rescaling $\pa$ by $\l$ will rescale one of the bases in the dual pair by $\l^{-1}$, and
thus will rescale $\mu$ by $\l^{-1}$. The pair $(\pa,\mu)$ will always satisfy $\pa(\mu) = k$, where $k \in \NM$ is the rank of the extension. Thus if $\pa$ and $\mu$ are known up to a
scalar, specifying either one will suffice to pin down the Frobenius extension structure.

Let us quickly note that for any Frobenius extension the trace map $\pa$ must be surjective.

Now we define the notion of a Frobenius hypercube, as in \cite{EWFrob}.

\begin{defn} Let $S$ be a finite set, and let $\ell$ be a function which assigns an integer to each subset of $S$. A \emph{Frobenius hypercube} is data of a (commutative graded) ring $R^I$
for each $I \subset S$, and whenever $I \subset J \subset S$, a Frobenius extension $R^J \subset R^I$ of degree $\ell(J)-\ell(I)$. The trace maps are required to satisfy $\pa^J_I \pa^K_J =
\pa^K_I$ whenever $K \subset J \subset I$.

A \emph{partial Frobenius hypercube} is the same definition, except restricted to subsets $I \subset S$ lying within a certain ideal (under the inclusion partial order). For example, one
might only consider proper subsets of $S$. \end{defn}

Our goal will be to define a partial Frobenius hypercube attached to a realization, though this will require slightly more data than the realization itself provides.

\subsection{Polynomials and invariants}
\label{subsec-invts}

\begin{assumption} \label{demazuresurjectivity} The maps $\a_s \co \hg \to \Bbbk$ and $\a_s^\vee \co \hg^* \to \Bbbk$ are surjective for each $s \in S$. This assumption is called
\emph{Demazure surjectivity}. \end{assumption}

If the row and column corresponding to $s$ in the Cartan matrix each contain an invertible element of $\Bbbk$, then Demazure surjectivity follows immediately. Thus Demazure surjectivity
holds for any exotic realization of $\sl_{n+1}$ for $n \ge 2$. Demazure surjectivity for exotic affine $\sl_2$ can be guaranteed by assuming that $2$ and $[2]$ generate the unit ideal in
$\Bbbk$ (and this is required when $\D$ and $\D^\vee$ are bases).

We let $R=\textrm{Sym}(\hg^*)$. When the simple roots form a basis for $\hg^*$, this is just a polynomial ring on the variables $\a_s$, $s \in S$. We give $R$ an even $\Z$-grading, so
that $\deg(\hg^*)=2$. There is an action of $W$ on $R$ by ring homomorphisms, extending the action on the linear terms.

Let $R^s$ denote the invariants of $R$ under $s \in S$. There exists an $R^s$-linear \emph{Demazure operator} $\pa_s \co R \to R^s$ of degree $-2$, defined on linear terms by $\pa_s =
\a_s^\vee \co \hg^* \to \Bbbk$, and satisfying the twisted Leibniz rule $\pa_s(fg) = \pa_s(f)g + s(f) \pa_s(g)$. It can be defined by \eqref{eq:partial}.

The Demazure operator $\pa_s$ equips the ring extension $R^s \subset R$ with the structure of a graded Frobenius extension (this relies on Demazure surjectivity). The invariant $\mu_s$ of
this extension structure is equal to $\a_s$. Thus while the rings $R$ and $R^s$ only depend on the action of $W$ on $\hg^*$, the Frobenius structure depends on the choice of simple roots
and coroots, and one should think of a simple realization as being the choice of Frobenius structures on $R^s \subset R$ for all $s \in S$.

We say that $I$ is \emph{finitary} when the parabolic subgroup $W_I \subset W$ is finite. Let $R^I \subset R$ denote the invariants under $W_I$. When $I$ is finitary, the Chevalley Theorem
states that $R$ is free over $R^I$ of finite rank (whereas when $W_I$ is infinite, $R^I$ has smaller transcendence degree, and the extension has infinite rank). Under some additional
assumptions, $R^I \subset R$ can be given the structure of a Frobenius extension of degree $\ell(I)$, where $\ell(I)$ is the length of the longest element $w_I \in W_I$.

\begin{claim} The operators $\pa_s$ satisfy the braid relations up to an invertible scalar. For any reduced expression $\un{w}_I = s_1 s_2 \ldots s_{\ell(I)}$ for $w_I$, one has a
composition of Demazure operators $\pa_{\un{w}_I} = \pa_{s_1} \cdots \pa_{s_d}$ of degree $-\ell(I)$, and this depends on the reduced expression only up to invertible scalar. Under some
additional assumptions on $\Bbbk$, any such map $\pa_{\un{w}_I}$ will equip $R^I \subset R$ with the structure of a Frobenius extension. Similarly, when $J \subset I$ are both finitary and
some additional assumptions hold, any reduced expression for $w_I w_J^{-1}$ will define an iterated Demazure operator giving a Frobenius trace $R^J \to R^I$. \end{claim}

\begin{remark} The statement that a map is a Frobenius trace is none other than the existence of dual bases, which requires the map to be surjective, and requires a certain determinant to
be invertible. What precisely this determinant is for the various traces $R^J \to R^I$ does not appear to be in the literature, so that we can not make the required assumptions as explicit
as Assumption \ref{demazuresurjectivity}. For the case of $\sl_3$, the assumptions can be guaranteed by assuming that $3$ is invertible (and this is required when $\D$ and $\D^\vee$ are
bases). \end{remark}

\begin{exercise} For the Frobenius trace $\pa^y_g$ defined in \S\ref{subsec-polynomials-sl3}, find dual bases for $R^y$ over $R^g$ under the assumption that $3$ is invertible.
\end{exercise}

\begin{defn} A \emph{Frobenius realization} is a simple realization equipped with a choice of Frobenius structure $\pa^J_I$ on $R^I \subset R^J$ for each finitary pair $J \subset I \subset
S$. (In particular, we assume that such a Frobenius structure exists.) The Frobenius trace $R \to R^s$ must simply be $\pa_s$, as determined by the simple realization. We require that
$\pa^J_I \pa^K_J = \pa^K_I$ whenever $K \subset J \subset I$. The Frobenius trace $\pa_I \co R \to R^I$ must be proportional to $\pa_{\un{w}_I}$ for some reduced expression of $w_I$.
\end{defn}

This provides a partial Frobenius hypercube. As the Frobenius trace maps $\pa^J_I$ are defined up to an invertible scalar, we can specify the structure by choosing the invariant $\mu^J_I$.
We do this by choosing a set of positive roots.

\begin{defn} A \emph{root realization} is a simple realization with the additional choice of a set of \emph{finite positive roots} $\Phi$. It must include the simple roots, and contain
exactly one vector collinear with each element of $W_I \cdot \D_I$ for $I$ finitary. In other words, whenever $\un{w} = s_1 \ldots s_d$ is a reduced expression within a finite parabolic
subgroup $W_I$, then $\Phi$ contains one element proportional to $s_1 s_2 \ldots s_{d-1} (\a_{s_d})$. Such an element is a \emph{positive root for $W_I$} and lives in the subset $\Phi_I
\subset \Phi$. \end{defn}

\begin{defn} Assume that $\pa_{\un{w}_I}$ is a Frobenius trace for some reduced expression of $w_I$, for each finitary $I$. Then any root realization will determine a unique Frobenius
realization for which the invariant $\mu_I$ of the extension $R^I \subset R$ is equal to the product of the positive roots for $W_I$. The invariant $\mu^J_I$ for $R^I \subset R^J$ is equal
to the product of the positive roots for $W_I$ which are not roots for $W_J$. We call this the \emph{associated Frobenius realization}. \end{defn}

Finally, let us note when these scalar ambiguities do not arise. The following is proven in \cite{EDihedral}.

\begin{defn} A simple realization is \emph{odd-balanced} if, whenever $m_{st}$ is odd for $s,t \in S$, then $a_{st} = a_{ts}$ satisfies the algebraic conditions to be equal to
$-(\z+\z^{-1})$ for a primitive $2m_{st}$-th root of unity $\z$. \end{defn}

\begin{claim} Whenever a realization is odd-balanced, there is a canonical choice of positive roots, and the Demazure operators satisfy the braid relations. \end{claim}

\subsection{Singular Bott-Samelson bimodules and diagrams}
\label{subsec-singBS}

Fix a Frobenius realization of $(W,S)$. One can now define $\SBSBim_{\textrm{alg}}$ and $\SSBim$ precisely as in \S\ref{subsec-summary}. One can also define a diagrammatic category
$\SBSBim_{\textrm{diag}}$ as in the previous sections (for $\sl_2$ and $\sl_3$).

Singular Soergel bimodules were defined by Williamson \cite{WilSSB}. When the realization is \emph{reflection faithful} (see \cite{WilSSB} or the introduction), the category $\SSBim$ is
well-behaved. It categorifies the Hecke algebroid, an ``idempotented" version of the Hecke algebra (and called the Schur algebroid in \cite{WilSSB}). The size of morphism spaces between
singular Soergel bimodules is governed by the standard trace on the Hecke algebroid; this is known as the \emph{Soergel-Williamson Hom formula}. When $\Bbbk=\C$ and the realization is
``geometric", the indecomposable objects in $\SSBim$ are (roughly speaking) the equivariant hypercohomologies of simple perverse sheaves on a Kac-Moody group, equivariant under various
parabolic subgroups, and thus they descend to the Kazhdan-Lusztig basis of the Hecke algebroid.

The situation is analogous to the discussion of $\Rep$ and $\Fund$ in the previous chapter. $\SSBim$ is the genuine object of interest, while $\SBSBim$ is the combinatorial replacement.
Finding a description of $\SBSBim$ by generators and relations has been elusive in general, and no progress has been made outside of rank $2$ and type $A$. Rank $2$ was accomplished in
\cite{EDihedral}. Williamson and I have a conjectural presentation of $\SBSBim$ in type $A$ which we hope will shortly be available. Thus the presentation $\SBSBim_{\textrm{diag}}$ for
affine $\sl_3$ is joint with Williamson; the proof of its correctness given in the appendix is original, but similar in style to arguments we attempted together.

When the realization is not reflection faithful, the category $\SSBim$ is not well-behaved, and will not categorify the Hecke algebroid. An example will be when $\Bbbk = \CM$ and $q \ne
\pm 1$ is a root of unity, where the realization ceases to be faithful. However, the diagrammatic category $\SBSBim_{\textrm{diag}}$ is still well-behaved, and its morphism spaces always
satisfy the Soergel-Williamson Hom formula. The presentation only encodes maps which occur generically; more morphisms may exist between algebraic Soergel bimodules when $q$ is a root of
unity, but these morphisms are not in our diagrammatic presentation.

\subsection{Maximally singular Bott-Samelson bimodules}
\label{subsec-restricted}

Let us quickly mention the particular sub-$2$-category of $\SBSBim$ that we will use in this paper. Fix a Dynkin diagram $\G$. Let $(W,S)$ be the affine Weyl group of $\G$, so that $S =
\tG$, and fix a Frobenius realization. Let $\Om \cong \trmv \subset S$ be as in \S\ref{subsec-om-versions}. Maximal finitary subgroups have the form $I=S \setminus s$ for $s \in S$, and
when $s \in \trmv$ the parabolic subgroup $W_I$ is isomorphic to the finite Weyl group.

\begin{defn} Let $\mSSBim$, the $2$-category of \emph{maximally singular Soergel bimodules}, denote the full sub-$2$-category of $\SSBim$ whose objects are subsets of the form $I = S
\setminus s$, for $s \in \trmv$. \end{defn}

We can now state the Soergel Satake theorem.

\begin{thm} (Soergel Satake) \label{thm:soergelsatake} There is a $2$-functor $\FC \co \RepOm(\gg^\vee) \to \mSSBim(\tG)$, which is a degree $0$ equivalence. \end{thm}

Henceforth in this chapter we assume $(W,S)$ is the affine Weyl group in type $A$. Then $\mSSBim$ has an even better combinatorial replacement than $\SBSBim \cap \mSSBim$.

\begin{defn} Let $\mSBSBim$, the $2$-category of \emph{maximally singular Bott-Samelson bimodules}, denote the full sub-$2$-category of $\SBSBim$ defined as follows. \begin{itemize} \item
The objects are subsets of the form $I = S \setminus s$, for $s \in \trmv$. \item The $1$-morphisms are monoidally generated by bimodules of the form ${}_{R^I} R^{I \cap J}_{R^J} (\ell(J)
- \ell(I \cap J))$, where $I = S \setminus s$ and $J = S \setminus t$ for $s \ne t \in \trmv$. \item The $2$-morphisms are all graded bimodule maps. \end{itemize} \end{defn}

\begin{claim} The Karoubi envelope of $\mSBSBim$ is $\mSSBim$. In other words, every summand of a singular Bott-Samelson bimodule between maximally singular parabolic subsets is a summand
of a maximally singular Bott-Samelson bimodule. \end{claim}

\begin{remark} This corresponds, under algebraic Satake, to the fact that all irreducible representations of $\sl_n$ are summands of tensor products of miniscule representations. For
general affine Weyl groups, more complicated bimodules are required to produce the non-miniscule fundamental representations. \end{remark}

Indecomposable singular Soergel bimodules are classified numerically, that is, by their behavior in the Grothendieck group. The proof of this claim is also numerical, showing that
maximally singular Bott-Samelsons also have a certain upper-triangularity property in the Grothendieck group. For the classification of indecomposable bimodules, see \cite{WilSSB}.

\begin{remark} \label{rmk:reverse} One might also call the objects of $\mSBSBim$ by the name \emph{reverse Bott-Samelson bimodules}. This is because, in their original context,
Bott-Samelson bimodules are $R$-bimodules obtained by restricting from $R$ to $R^s$ and inducing back to $R$. In other words, we begin at the minimal finitary parabolic subset $\emptyset$,
move to a subminimal one $\{s\}$, and come back to the minimal one. For reverse Bott-Samelson bimodules we begin at a maximal finitary parabolic subset $S \setminus s$, move to a
submaximal one $S \setminus \{s,t\}$, and return to a (different) maximal one $S \setminus t$. \end{remark}

We can now state the algebraic Satake theorem. We now reintroduce the subscript $q$ to indicate that we are taking an exotic affine $\sl_n$ realization.

\begin{thm} (Algebraic Satake) \label{thm:algebraicsatake} There is a $2$-functor $\FC \co \FundOm_q \to \mSBSBim_q$, which is a degree $0$ equivalence. \end{thm}

\section{Reformulating Geometric Satake}
\label{sec-reformulate}

The reformulation from geometric Satake to Soergel Satake takes the form of three ``replacements:" \begin{itemize} \item Replacing monoidal categories with $2$-categories, \item Replacing
the affine Grassmannian with a partial flag variety for the affine Kac-Moody group, \item Replacing perverse sheaves with their global sections (i.e. their equivariant hypercohomology).
\end{itemize} In the sections to come we will present this flow of ideas for a general complex semisimple lie algebra $\gg$. Nothing used is new or should be unfamiliar to the experts.

The transformation from Soergel Satake to algebraic Satake consists of replacing each side of the equivalence with a combinatorial additive subcategory. The technology for this procedure
only exists currently in type $A$. This final replacement is discussed in greater detail in sections \S\ref{subsec-webs-philosophy} and \S\ref{subsec-restricted}, but we also mention it
briefly here.

\subsection{Geometric Satake}
\label{subsec-gsatakebackground}

We begin by recalling the results which are packaged under the name ``geometric Satake." We refer the reader to \cite{GinzburgGS} for better and more thorough introduction to this topic.
We assume for the moment that all categories are $\C$-linear. For a group $H$ acting on a space $X$, we let $\PC_H(X)$ denote the category of $H$-equivariant perverse sheaves on $X$. When
two groups $H$ and $K$ act on $X$, on the left and right respectively, we write $\PC_{(H,K)}(X)$ for the category of biequivariant perverse sheaves.

Fix $T\subset G$ a maximal torus in a simple algebraic group, with Lie algebras $\hg \subset \gg$. Let $G^\vee$ (resp. $\gg^\vee$)
denote the Langlands dual group (resp. Lie algebra). Let $\KC=\C((t))$ and $\OC=\C[[t]]$. As it is usually stated, geometric Satake is an equivalence between $G^\vee(\C)$-rep and the
category of $G(\OC)$-equivariant perverse sheaves on the affine Grassmannian $\Gr_{G} = G(\KC)/G(\OC)$.

This equivalence intertwines the (non-equivariant) hypercohomology functor with the forgetful functor: the hypercohomology of an irreducible perverse sheaf is isomorphic, as a vector
space, to the underlying vector space of the corresponding irreducible $G^\vee$ representation. In this context, the hypercohomology and forgetful functors are called \emph{fiber
functors}. The cohomological grading on hypercohomology transforms to the weight grading on the representation (with respect to some naturally-defined regular semisimple element). More
interesting features arise when one considers the $T$-equivariant hypercohomology: instead of a vector space, one gets a graded module over the $T$-equivariant cohomology of a point, which
is the polynomial ring $\C[\hg^*]$ (graded with $\deg \hg^* = 2$). Taking the quotient by the augmentation ideal (i.e. killing positive degree polynomials), one obtains the usual
hypercohomology, though other specializations also hold interest. The usual tensor product on representations is intertwined with a convolution functor on perverse sheaves, though the
geometric proof that convolution is symmetric monoidal is quite technical.

\begin{remark} \label{rmk:caveat2} This long remark continues the discussion from the introduction (Remark \ref{rmk:caveat}) about the validity of discussing geometric
Satake without the symmetric structure.

Given two semisimple monoidal categories with the same classification of indecomposables and the same decomposition of tensor products (i.e. the same Grothendieck ring), when are they
equivalent? The answer is encoded in the associativity isomorphism. If these monoidal categories are also symmetric, then checking equivalence of symmetric monoidal categories also
requires checking the commutation isomorphism. See \cite[Problem 1.42.8 and following]{EtingofTensorNotes} for more details.

The numerical underpinnings of geometric Satake were shown by Lusztig in \cite{LusztigGS}. The characters of the irreducible perverse sheaves are encoded in Kazhdan-Lusztig polynomials for
the affine Hecke algebroid. According to geometric Satake, these characters (evaluated at $v = 1$) should agree with the weight multiplicities in the corresponding irreducible
$G^\vee$-representation. This was shown by Lusztig, who also demonstrated the link with orbit closures in the affine Grassmannian. Moreover, these results already encode a great deal about
the monoidal structure on perverse sheaves. When one convolves two irreducible perverse sheaves, the BBD Decomposition Theorem \cite{BBD} implies that the result splits into irreducibles
according to its character, just as the splitting of a tensor product of $G$-representations is determined by its weight multiplicities. Thus perverse sheaves form a semisimple monoidal
category with the same Grothendieck ring as $G^\vee(\C)$-representations. The fiber functor (implicit also in Lusztig's work) also commutes with the convolution structure, making it easy
to check the associativity isomorphism (note: this argument relies on semisimplicity). This is done explicitly in Ginzburg's preprint \cite{GinzburgGS}. With Lusztig's results and the
Decomposition Theorem in hand, the only further content in Geometric Satake is the comparison of commutativity isomorphisms. Without this comparison one has only an equivalence of
(vanilla) monoidal categories; this is what we accomplish in this paper.

An elementary proof of the commutativity isomorphism was attempted by Ginzburg in \cite{GinzburgGS}, but is not believed to be correct. A correct (and non-elementary) proof was given by
Mirkovic and Vilonen \cite{MirkVil}, following ideas of Drinfeld. These results work with a ``global approach," using the fusion product and the Beilinson-Drinfeld Grassmannian. This is
the style also used by Gaitsgory in his quantum geometric Satake equivalence \cite{Gaits}. The author does not know how the global approach connects to the Soergel picture.

We do rely on Lusztig's numerical results, though we do not use the Decomposition theorem. Instead, we compare two monoidal categories by looking inside at strictly-monoidal additive
subcategories (where the associativity isomorphism is by definition the identity map). To show two such categories are equivalent one has no convenient tricks (to the author's knowledge),
but must compare the morphism algebras in full, as we did above. In essence, we have converted the questions of decomposition and associativity into a question about morphisms in this
subcategory. The proof we use does not use the Decomposition Theorem or semisimplicity in any way, relying on Soergel's categorification results but not on the Soergel conjecture. In other
words, we never need to study indecomposable singular Soergel bimodules and their numerical properties, only those bimodules which correspond to tensor products of fundamental
representations. This is one advantage of working directly with generators and relations.

A proof independent of the Decomposition Theorem is required for the quantum algebraic Satake equivalence (and not just because there is no corresponding geometry)! After all, when $q$ is
a root of unity, the analog of the Soergel conjecture will fail, and the indecomposable objects will have the ``wrong" characters, just as for representations of the (Lusztig form of the)
quantum group at a root of unity. Nevertheless, tensor products of fundamental representations (and their Soergel analogs) continue to have the same numerics. \end{remark}

\subsection{From monoidal categories to 2-categories}
\label{subsec-reform-monoidalto2}

Our first step will be to apply a sequence of tautologies in order to reformulate geometric Satake as an equivalence of (strict) abelian $2$-categories. In \S\ref{subsec-om-monoidal} it is
explained how to begin with representations of $\gg^\vee$ and obtain $\RepOm$, an $\Om$-$2$-category. This will be one side of our equivalence. By conventional geometric Satake, the
geometric category paired with representations of $\gg^\vee$ would be perverse sheaves on $\Gr_{G_{\adj}}$ in adjoint type. So we begin by extracting a $2$-category from
$G_{\adj}$, after which we will explain how the same data is encoded within the affine Grassmannian for any $G$.

Suppose that $G = G_\adj$ has adjoint type, so that $\pi_0(G(\KC))$ is the group $\Om$ defined in \S\ref{subsec-om-versions}. As $G(\OC)$ is connected, this is also the component group of
$\Gr_{G}$. Let $\Gr_0$ denote the \emph{nulcomponent}, i.e. the component associated to the identity in $\Om$. Since $\Om$ is abelian, conjugation by any element will preserve the
nulcomponent $G(\KC)_0$. Any irreducible perverse sheaf is supported on a single component, and perverse sheaves on different components admit no morphisms or extensions. Convolution will
act on the ``component support" of a perverse sheaf precisely via multiplication in $\Om$, so that $\PC_{G(\OC)}(\Gr_{G})$ already forms an $\Om$-graded monoidal category, from which one
can construct an $\Om$-$2$-category.

We choose to think of this $\Om$-$2$-category in a different way, altering equivariant structure rather than working on different components. Let us pick an equivariant perverse sheaf
$\FC$ supported on component $\xi \in \Om$. Let $x \in G(\KC)$ be an element in the component $\xi$. Then $x^* \FC$ is supported on the nulcomponent $\Gr_0$, and the
$G(\OC)$-equivariant structure on $\FC$ produces a natural $(x^{-1} G(\OC) x)$-equivariant structure on $x^* \FC$. The group $(x^{-1} G(\OC) x)$ depends on $\xi$, not on the
choice of $x$, so we call this group $G_{\xi}$. Therefore, instead of considering $\PC_{G(\OC)}(\Gr_{G})$, we can consider $\bigoplus \PC_{G_\xi}(\Gr_0)$ for all
$\xi \in \Om$.

\begin{example} \label{ex:sl4groups} Suppose that $\gg = \sl_4$, so that $G_\adj = \PM GL_4$, $G_\sc = SL_4$, and $\Om = \ZM/4\ZM$. The four subgroups $G_\xi \subset G_{\adj}(\KC)$ are as follows:
\[G_{\overline{0}} = \left( \begin{array}{cccc} \OC & \OC & \OC & \OC \\ \OC & \OC & \OC & \OC \\ \OC & \OC & \OC & \OC \\ \OC & \OC & \OC & \OC \end{array} \right),
\qquad
G_{\overline{1}} = \left( \begin{array}{cccc} \OC & t^{-1} \OC & t^{-1} \OC & t^{-1} \OC \\ t \OC & \OC & \OC & \OC \\ t \OC & \OC & \OC & \OC \\ t \OC & \OC & \OC & \OC \end{array}
\right), \]
\[G_{\overline{2}} = \left( \begin{array}{cccc} \OC & \OC & t^{-1} \OC & t^{-1} \OC \\ \OC & \OC & t^{-1} \OC & t^{-1} \OC \\ t \OC & t \OC & \OC & \OC \\ t \OC & t \OC &
\OC & \OC \end{array} \right),
\qquad
G_{\overline{3}} = \left( \begin{array}{cccc} \OC & \OC & \OC & t^{-1} \OC \\ \OC & \OC & \OC & t^{-1} \OC \\ \OC & \OC & \OC & t^{-1} \OC \\ t \OC &
t \OC & t \OC & \OC \end{array} \right). \]
These groups are obtained from $G_{\overline{0}}$ by conjugation by a diagonal matrix with entries $(t,\ldots,t,1,\ldots,1)$. This diagonal matrix is always an element of $\PM GL_4
(\KC)$, though not usually an element of $SL_4(\KC)$. Note that the intersection $\cap_{\xi} G_\xi$ is an Iwahori subgroup, a fact which is special to type $A$. \end{example}

Perverse sheaves which are $G(\OC)$-equivariant on $G(\KC)/G(\OC)$ have a natural monoidal structure from convolution. The collection of $G_\xi$-equivariant perverse
sheaves on the nulcomponent of $G(\KC)/G_{\overline{0}}$ does not have an obvious monoidal structure (it is a module over the previous monoidal category). However, by considering
the collection of all $G_\xi$-equivariant perverse sheaves on the nulcomponent of $G(\KC)/G_{\eta}$ for $\xi, \eta \in \Om$, one obtains a $2$-category where convolution
again makes sense. We denote this $2$-category by $\PPP$.

For mental simplification, we choose to work with $\PC_{(G_\xi,G_\eta)}(G(\KC)_0)$, i.e. bi-equivariant sheaves on the nulcomponent of $G(\KC)$, instead of $\PC_{G_\xi}(G(\KC)_0/G_\eta)$,
i.e. equivariant sheaves on the nulcomponent of the quotient. Because $G_\xi$ orbits on $G(\KC)$ are infinite-dimensional, it is difficult to comprehend what is meant by equivariant
perverse sheaves on $G(\KC)$. However, in any sense in which such things are defined, there will be an equivalence between these two categories. We choose the former mostly for the
notational symmetry, and because it makes the convolution structure ``seem" more natural.

Now let us explain how the same $2$-category $\PPP$ can be extracted from any Lie group $G$ lifting $\gg$. The natural map $G \to G_{\adj}$ is a Galois cover with
kernel $Z$, the center of $G$. It induces a map $G(\OC) \to G_{\adj}(\OC)$ whose kernel is also $Z$, and a map $G(\KC) \to G_{\adj}(\KC)$ whose kernel, restricted
to the nulcomponent, is also $Z$. Thus the induced map $\Gr_{G} \to \Gr_{G_\adj}$ is an isomorphism on the respective nulcomponents. Thus the space on which we consider our
sheaves does not depend on the choice of group.

The group $\Om$ acts on $G(\KC)_0$ by group automorphisms, for any $G$. One can still think of this action as ``conjugation" by an element in $G_{\adj}(\KC)$, but this
automorphism may be an outer automorphism for $G(\KC)$. Thus one still has groups $G_{\xi} = \xi(G(\OC))$ for each $\xi \in \Om$. These groups will also be the preimages of
the corresponding groups $G_{\xi} \subset G_{\adj}(\KC)_0$, and their orbits on $G(\KC)$ will be the preimages of orbits on $G_{\adj}(\KC)$.

\begin{example} Continuing example \ref{ex:sl4groups} above, one can define $G_{\overline{k}}$ exactly as above for each $0 \le k \le 3$, with the assumption that one only considers
matrices of that form lying within $G$. Conjugation by the diagonal matrix with entries $(t,\ldots,t,1,\ldots,1)$ is a (potentially outer) automorphism of $G(\KC)$. \end{example}

Now compare the equivariant categories $\PC_{G_\xi}(G(\KC)_0/G_\eta)$. As the group changes, the space $G(\KC)_0/G_\eta$ is unchanged, but the group which acts
$G_\xi$ will be altered by some central kernel $Z$. Clearly $Z$ acts trivially, being contained in both $G_\xi$ and $G_\eta$. However, $Z$ is in the nulcomponent of any
stabilizer, meaning that it will also act trivially on any equivariant structure. Thus the equivariant category is also independent of the choice of $G$.

Taking the sum over all $\PC_{G_\xi}(G(\KC)_0/G_\eta)$, one obtains a $2$-category which is equivalent to $\PPP$. This $2$-category is larger than merely $\PC_{G(\OC)}(\Gr_{G})$; those
$G_\xi$ which are genuinely conjugate to $G(\OC)$ will encode perverse sheaves on other connected components, while the ones which are only outer-conjugate will allow for the extra data.

\subsection{From affine Grassmannians to Kac-Moody groups}
\label{subsec-reform-grasstoKM}

Consider the Kac-Moody group $G_\aff$ associated to the affine Dynkin diagram $\tG$, as in \cite{Kumar:book}. This is an analogue of $G(\KC)$. The essential difference is that $G_\aff$ has
a torus whose dimension is bigger by two: it has a copy of $\C^*$ acting on $G(\KC)$ by loop rotation (sending $t \mapsto \l t$), and a copy of $\C^*$ arising from a central extension.

The set of removable vertices $\trmv$ was defined in \S\ref{subsec-om-versions}; recall that it is a $\Om$-torsor. For any subset of $\tG$ one has a corresponding parabolic subgroup of
$G_\aff$, so for each $s \in \trmv$ one has a parabolic subgroup $P_s$ corresponding to $\tG \setminus s$. Because each of these parabolic subgroups also contains the affine torus, the
extra dimensions cancel in the quotient space $G_\aff / P_s$, which is isomorphic to $\Gr^\vee_0$ for any $s \in \trmv$. However, the structure of $P_t$-equivariance on $G_\aff/ P_s$ is
stronger than the structure of $G^\vee_\xi$-equivariance, because of the larger torus. The central extension $\C^*$, being central, will act trivially on $G_\aff / P_s$, and being
connected, will therefore act trivially on any equivariant structure. However, the loop rotation $\C^*$ will add new structure.

Thus our new $2$-category of interest is the collection of $P_t$-equivariant perverse sheaves on $G_\aff/P_s$ for various $s,t \in \trmv$, or alternatively, $\PC_{(P_t,P_s)}(G_\aff)$. This
is a souped-up version of the original $2$-category, having an extra ``dimension" of equivariant structure.

Note that in type $A$, because $\trmv = \tG$, the intersection of all the $P_s$ will be precisely the Borel subgroup. This is analogous to the situation in Example \ref{ex:sl4groups},
where the intersection of the $G_\xi$ was an Iwahori subgroup.

\subsection{From perverse sheaves to bimodules}
\label{subsec-reform-globalsections}

We now describe some results of Soergel and Williamson, in order to motivate the theorem of H\"arterich we shall use.

Recall that the $T$-equivariant cohomology of a point $H_T(\pt)$ is naturally isomorphic to $R_\fin = \CM[\hg_\fin^*]$, a polynomial ring generated by dual Cartan subalgebra. We
temporarily denote this ring $A$, not to be confused with the ring $R$ we have used elsewhere in this paper. For $T$-spaces with a cellular filtration, passing from the $T$-equivariant
cohomology $H_T(X)$ to the ordinary cohomology $H(X)$ is obtained by tensoring over $H_T(\pt)$ with the one-dimensional module $A / A_+ = \CM$, by the localization theorem. There is a
natural action of the finite Weyl group $W_\fin$ on $A$, and the $G$-equivariant cohomology of a point is $A^{W_\fin}$. For the parabolic subgroup $P_I$ associated to $I \subset S_\fin$,
the $P_I$-equivariant cohomology of a point is $A^I = A^{W_I}$.

Therefore, given a $(P_I,P_J)$-equivariant perverse sheaf on $G$, its equivariant hypercohomology (i.e. equivariant pushforward to a point) will naturally be a graded $(A^I,A^J)$-bimodule.

\begin{thm} The functor of equivariant hypercohomology from $\PC_{(P_I,P_J)}(G)$ to graded $(A^I,A^J)$-bimodules is fully faithful on semisimple objects, for any $I,J \subset S$. (Note
that this functor sends higher extensions to bimodule maps of nonzero degree, as the rings $A^I$ are graded homologically.) The images of the semisimple objects are known as \emph{singular
Soergel bimodules} for $W_\fin$. \end{thm}

\begin{remark} Unfortunately, this precise version of the theorem does not appear in the literature, so we will briefly discuss the history and sketch the proof. Soergel first studied the
case $P_I = P_J = B$. A proof using ordinary (rather than equivariant) hypercohomology is in \cite{Soe1}, while its equivariant analog is in \cite[Proposition 3.4.4]{Soe2001}. The general
idea of these proofs is the same. First, one shows that pullbacks and pushforwards between partial flag varieties (or induction/restriction of equivariance) is sent by hypercohomology to
induction and restriction of bimodules. This implies that the images of semisimple objects are Soergel bimodules. Next, one uses a parity-vanishing argument to deduce that the spectral
sequence computing morphisms between hypercohomologies actually degenerates on the first page, implying that hypercohomology is faithful. Finally, one uses the Soergel Hom formula to
deduce that the dimensions of morphism spaces agree. For the general case a similar argument will apply, using Williamson's generalization of the Soergel Hom formula \cite{WilSSB}; though
not in the literature, this proof is known to experts. \end{remark}

The theorem was upgraded to affine Kac-Moody groups by Soergel's student H\"arterich \cite{Harterich}. Now one works with the Kac-Moody group $G_\aff$ and its parabolic subgroups $P_I$
for $I \subset S$.

\begin{thm} The functor of equivariant hypercohomology from $\PC_{(P_I,P_J)}(G_\aff)$ to graded $(H_{P_I}(\pt),H_{P_J}(\pt))$-bimodules is fully faithful on semisimple objects,
for any $I,J \subset S$ with $W_I, W_J$ finite. The images of the semisimple objects are known as \emph{singular Soergel bimodules} for $W_\aff$. \end{thm}

Let $T_\aff$ denote the torus of $G_\aff$, let $T'$ denote the finite torus extended by loop rotation, and let $S$ denote the set of affine reflections. Then $H_{T'}(\pt)$ is naturally
isomorphic to the polynomial ring $R=\CM[\hg_\aff^*]$. This is the ring $R$ we have used throughout this paper, defined in Example \ref{ex:stdaffine}, or in
\S\ref{subsec-exoticrealization} at $q=1$. On the other hand $R_\aff = H_{T_\aff}(\pt)$ is a larger ring, corresponding to the self-dual $n+2$-dimensional realization also mentioned in
Example \ref{ex:stdaffine}. The subrings $H_{P_I}(\pt)$ are isomorphic to $R_\aff^I$, so that $H_{P_s}(\pt)$ is isomorphic to $R_\aff^{W_\fin}$ for any $s \in \trmv$.

The difference between $R_\aff$ and $R$ is a polynomial ring in one variable, on which $W$ acts trivially: namely, $H_{\C^*}(\pt)$ for the central extension. All computations done in this
paper work for $R_\aff$ as well as for $R$. Conversely, the theorem above will apply equally well to bimodules over the subrings of $H_{P_I}(\pt)$ which ignore the central extension.
Henceforth, we will work with singular Soergel bimodules for $R$ rather than for $R_\aff$.

The singular Soergel bimodules defined in \S\ref{subsec-summary} agree with the singular Soergel bimodules in this theorem. Therefore, the (semisimple part of the) $2$-category of
interest, $\PC_{(P_t,P_s)}(G_\aff)$, is equivalent to $\mSSBim$, defined in \S\ref{subsec-restricted}. In fact, every perverse sheaf in this $2$-category happens to be semisimple (warning:
higher extensions do exist).

Here we see the payoff of having replaced the affine Grassmannian with the Kac-Moody group. The $G(\OC)$-equivariant cohomology of the point is the same as the $G$-equivariant cohomology,
which is $R_\fin^{W_\fin}$. Meanwhile, the $P_s$-equivariant cohomology of the point is $R_\aff^{W_\fin}$, a larger ring. One could not repeat the constructions of the previous chapters
with $R_\fin$, because only one copy of $W_\fin$ acts on it, while for $R_\aff$ or $R$ there are a number of copies of $W_\fin$ inside $W_\aff$ which can act.

\subsection{Combinatorial subcategories}
\label{subsec-combointro}

When investigating an additive (monoidal) category, it may be difficult to compute all the morphisms between indecomposable objects. However, there may be a class of objects (closed under
the tensor product) for which the computation of morphisms is tractable, and for which every indecomposable object appears as a summand. We call the corresponding full subcategory a
\emph{combinatorial replacement} for the original additive category.

A discussion of this philosophy and the known results for the special case of $\Rep(\gg)$ or $\RepOm(\gg)$ can be found in \S\ref{subsec-webs-philosophy}. In that section, it is discussed
that tensor products of fundamental representations form a combinatorial replacement in type $A$, where the morphisms between such tensor products can be described by $\sl_n$-webs.
Unfortunately, morphisms between tensor products of fundamental representations remain undiscovered outside of type $A$ and rank $2$. The description of webs for $\sl_2$ and $\sl_3$ is
given explicitly in \S\ref{subsec-webs-sl2} and \S\ref{subsec-webs-sl3} respectively.

The same philosophy applies to perverse sheaves, and in fact was a strong motivation for Soergel's original definition of Soergel bimodules. Soergel originally worked with the finite flag
variety, or $\PC_{(B,B)}(G)$. Finding irreducible perverse sheaves is difficult, but the decomposition theorem implies that a given irreducible perverse sheaf is a summand of the
pushforward of the constant sheaf from a proper resolution of singularities. Schubert varieties in the finite flag variety have ``combinatorial" resolutions known as \emph{Bott-Samelson
resolutions}, and the corresponding pushforwards are known as \emph{Bott-Samelson sheaves}. Continuing to pushforward these sheaves to a point, one obtains the Bott-Samelson bimodules
defined in Remark \ref{rmk:reverse} or in \cite{Soe5}. Bott-Samelson bimodules are easy to compute with, and thus form an excellent combinatorial replacement for perverse sheaves. As
the philosophy in \S\ref{subsec-webs-philosophy} points out, returning from an algebraic description of Bott-Samelson bimodules to a similar description of Soergel bimodules is very
difficult, and amounts to a computation of certain idempotents.

Williamson \cite{WilSSB} continued to treat the case of partial flag varieties (or rather, equivariance under various parabolic subgroups), and defined singular Soergel bimodules. Singular
Bott-Samelson bimodules, defined in \S\ref{subsec-summary} in terms of induction and restriction bimodules, form a combinatorial replacement in this context. For the specific case of
$\PC_{(P_s,P_t)}(G_\aff)$ in type $A$, it turns out that there is also a further replacement which works, discussed in \S\ref{subsec-restricted}.

\appendix
\section{Diagrammatic proofs for $\sl_3$}
\label{sec-proofs}

In this appendix, our primary goal is to prove the following theorem, which uses terminology from \S\ref{sec-ssbim}, and recalls the discussion of \S\ref{subsec-singBS}.

\begin{thm} \label{superthm-sl3} Consider a Frobenius realization of affine $\sl_3$ with exotic Cartan matrix. Let $\SBSBim_{\textrm{diag}}$ denote the $2$-category defined in
\S\ref{subsec-ssbimdiag-sl3}, and let $\SBSBim_{\textrm{alg}}$ denote the $2$-category defined in \S\ref{subsec-summary}. Then $\SBSBim_{\textrm{diag}}$ categorifies the Hecke algebroid,
and satisfies the Soergel-Williamson Hom formula. There is a faithful, essentially surjective $2$-functor $\GC \co \SBSBim_{\textrm{diag}} \to \SBSBim_{\textrm{alg}}$ (defined at the end
of \S\ref{subsec-ssbimdiag-sl3}). If the realization is reflection faithful, then $\GC$ is an equivalence. \end{thm}

Recall that when the realization is reflection faithful, Williamson's results imply that $\SBSBim_{\textrm{alg}}$ will also categorify the Hecke algebroid, and satisfy the
Soergel-Williamson Hom formula. Therefore, faithfulness of $\GC$ will imply fullness. When the realization is not even faithful (e.g. when $q$ is a root of unity in $\Bbbk$),
$\SBSBim_{\textrm{alg}}$ will fail to categorify the Hecke algebroid, and $\GC$ will not be full. Theorem \ref{superthm-sl3} clearly implies Claim \ref{claim:sl3ssbim}, and the reflection
faithful condition explains the restriction on $\Bbbk$ given in that claim.

Our second goal is to prove Theorem \ref{thm:sl3main}, which states that the functor $\FC \co \FundOmq \to \mSBSBim_q$ is a degree zero equivalence. This is independent of the choice of
realization (e.g. it still holds when $q$ is a root of unity).

Both results were proven for $\sl_2$ in \cite{EDihedral}. After setting up and proving some diagrammatic lemmata, the main arc of the proof is completely analogous to the $\sl_2$ case.
Like all diagrammatic proofs, it is extremely delicate.

Henceforth, to simplify notation, $\DC$ will denote the diagrammatic $2$-category $\SBSBim_{\textrm{diag}}$, and $\CC$ will denote the algebraic category $\SBSBim_{\textrm{alg}}$. A
\emph{diagram} refers to a singular Soergel diagram.

\subsection{Preliminaries}
\label{subsec-preliminaries}

We recall and amplify our conventions from \S\ref{sec-sl3-diag}. Let $W$ be the affine Weyl group of $\sl_3$. Let $S = \{r,b,y\}$ denote the set of \emph{primary colors}, which can be
combined to form the \emph{secondary colors} $\{g,o,p\}$. The color \emph{brown}, which combines all three primary colors, is forbidden as a region label in a diagram. For a subset $K
\subset S$, a diagram is said to only have the primary colors in $K$ if every region label is contained in $K$.

\begin{defn} \label{defn:DCK} Let $K \subsetneq S$ be a subset of the primary colors. Let $\DC(K)$ denote the $2$-category defined as in Definition \ref{defn:SBSBim-sl3} except that
\begin{itemize} \item $2$-morphisms may only have primary colors in $K$ (this also restricts the objects and $1$-morphisms), \item one only imposes the relations with primary colors in
$K$. \end{itemize} Let $\CC(K)$ denote the algebraic subcategory of singular Bott-Samelson bimodules defined as in Definition \ref{defn:ssbimtrue}, except that one only allows sets $I,J$
contained in $K$. \end{defn}

Note that both categories use the same base ring $R$, defined in \S\ref{subsec-polynomials-sl3}, regardless of the subset $K$ (we consider polynomials to be colorless).

There is a natural $2$-functor $\iota_K \co \DC(K) \to \DC$, which we will eventually show is fully faithful. At the moment, it is not even clear that $\iota_K$ is faithful, because $\DC$
has more relations. We prove this very soon - that relations between diagrams with extraneous primary colors do not affect diagrams without those primary colors. It is far less obvious
that $\iota_K$ should be full - that diagrams without certain primary colors on the boundary should be in the span of diagrams without those primary colors anywhere.

The Dihedral Cathedral \cite{EDihedral} is an in-depth study of Soergel bimodules for dihedral groups. Its appendix treats odd-unbalanced realizations such as the exotic realization of
affine $\sl_3$. We will quote this paper for many ``two-color" results in this proof, such as the following proposition.

\begin{prop} \label{prop:2colorGC} Let $K \subsetneq S$. Then the analogous $2$-functor $\GC_K$ is well-defined, full, and faithful from $\DC(K)$ to $\CC(K)$. \end{prop}

\begin{proof} Observe that any exotic realization of affine $\sl_3$ restricts to a reflection faithful representation of any dihedral parabolic subgroup. Now the result follows from the
Cathedral. \end{proof}

All the relations which define $\DC$ involve at most two primary colors. Therefore, Proposition \ref{prop:2colorGC} implies that these relations hold algebraically in $\CC$, and the
$2$-functor $\GC$ is well-defined. Moreover, $\GC_K = \GC \circ \iota_K$ is faithful, which implies that $\iota_K$ must be faithful.

The proof that $\iota_K$ is full will comprise the bulk of this appendix, after which certain categorification considerations will complete the proof of Theorem \ref{superthm-sl3}. It
amounts to a ``color-removal algorithm."

We say that a diagram \emph{reduces} to another class of diagrams if it can be rewritten, using the relations, as a linear combination of diagrams in this class. A region is
\emph{external} if it touches the boundary, and \emph{internal} otherwise. A diagram has \emph{featureless boundary} if it has a single exterior region. A \emph{loop} is a closed
$1$-manifold in a diagram (which may intersect strands of other colors), and may be oriented either clockwise or anticlockwise.

\begin{prop} \label{prop:colorremoval} Suppose that $K \subsetneq S$. A diagram whose external regions only contain primary colors in $K$ will reduce to diagrams where every region only
contains primary colors in $K$. \end{prop}

This is a restatement of the fullness of $\iota_K$. We now summarize the argument.

Suppose that $K = \{b,y\}$. The region colors we wish to remove, namely $r, o, p$, we will call ``reddish." Recall that we have also colored the strands in our diagrams, so that a red strand
separates reddish regions from non-reddish regions. Any diagram without external reddish regions but with some internal reddish region must have a clockwise red loop. Color-removal is
equivalent to the statement that diagrams without clockwise loops will span all diagrams in $\DC$.

A diagram without any crossings or closed $1$-manifolds is called a \emph{crossingless matching (with boxes)}. It is built entirely of cups, caps, and boxes. A crossingless matching with
featureless boundary must be a box. One nice feature of crossingless matchings on simply-connected domains is that they necessarily have a cup/cap on the boundary (i.e. two adjacent points
connected by a strand). However, this does not hold true of crossingless matchings on non-simply-connected domains, such as a configuration of radii on an annulus.

\begin{lemma} \label{lem:allreddish} Suppose that every region in a diagram (on a simply-connected domain) is reddish. Then the diagram reduces to a crossingless matching. \end{lemma}

\begin{proof} As the primary color red is never removed there can be no red strands, only blue and yellow strands. These can never cross, lest there exist a brown region. So we need only
show that any closed loop can be removed. Suppose there is a closed blue or yellow loop. Choosing an innermost loop (i.e. no loops in its interior; this makes sense since loops can not
cross), we replace it and its interior with a box using \eqref{bigdemazure} or \eqref{bigcircle}. Repeating this process, we remove all closed components. \end{proof}

To remove a clockwise red loop, one step of the algorithm will be to simplify its interior. A naive approach would go as follows. Suppose one could remove all closed loops from the
interior of the red loop. Then this interior is a diagram where every region is reddish, and by the above lemma this interior is a crossingless matching and must have a cup/cap on the
boundary. Applying \eqref{R2moves} one should be able to pull this cup/cap out of the yellow loop, reducing the number of strands which intersect the red loop, and allowing an inductive
argument. (Boxes are easy to deal with and do not complicate the argument overmuch, thanks to \eqref{singR2nonoriented1var} and similar relations.)

Although this style of argument is useful, the supposition is too naive. A complicating feature of the interior is that it may contain numerous anticlockwise red loops, and anticlockwise
loops cannot be removed in general! The resulting reddish domain is no longer simply connected, and we can not proceed as above. Instead, the interior of an anticlockwise loop can be
simplified, and this is where we must begin a more sophisticated argument. For this purpose, we prove some lemmata about Bott-Samelson objects.

\subsection{Bott-Samelson objects}
\label{subsec-bottsam}

\begin{defn} A \emph{BS object} in $\SBSBim$ is an object of the form $\mt s_1 \mt s_2 \mt \cdots \mt s_d \mt$. We say that a diagram on the disk has \emph{BS
boundary} if the boundary (read around the circle) is a BS object. \end{defn}

The following map has an alternating two-color BS boundary, and will be denoted $v_k$. \begin{equation} {
\labellist
\small\hair 2pt
 \pinlabel {$v_k$} [ ] at 128 36
\endlabellist
\centering
\ig{1}{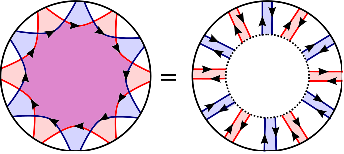}
}\end{equation} Pictured is the case $k=6$. When $k=0$, it looks like below. \igc{1}{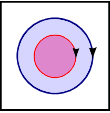} In $v_0$, switching the roles of blue and red would make no difference, as \eqref{bigcircle}, \eqref{slidepolyintoblue}, and \eqref{circleis3} imply that $v_0$ is equal to the polynomial $\mu_p$. We continue to denote the map $v_k$, regardless of which two primary colors are used.

The diagram $v_3$ is special, which relates to the fact that the order of $rb$ in $W$ is $3$. In \cite[section 6.1.2]{EDihedral} it was shown that every $v_k$ for $k \ge 4$ can be
rewritten as a composition of the maps $v_3$, together with cups and caps. It was also shown that every $v_k$ for $k < 3$ can be rewritten in terms of crossingless matchings. This can be
seen directly: for $k=2$ this is \eqref{squaresbsb}, for $k=1$ this is \eqref{singR2nonoriented2}, and for $k=0$ this follows as above.

The following claim is one of the results in the Cathedral.

\begin{claim} \label{BSmaps} Any diagram with at most two primary colors and BS boundary reduces to a composition of cups, caps, boxes, and the diagram $v_3$. If the BS boundary
is empty (i.e. featureless and white), then the diagram reduces to a box. If the BS boundary is non-empty then after reduction the diagram will either have a cap/cup on the
boundary, or will be equal to $v_k$ for some $k \ge 3$ (with boxes allowed in any white regions). \igc{1}{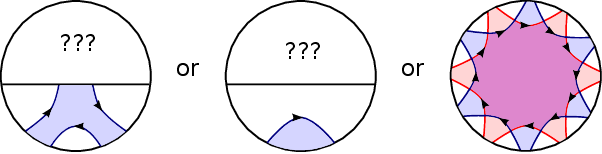} \end{claim}

In particular, if the BS boundary does not alternate between the two colors, but instead has repetition $\mt b \mt b \mt$, then there is guaranteed to be a cup/cap on the boundary as in
the first option above. One way to see this is to apply \eqref{circforcesamealt} to the intermediate $b \mt b$.

\begin{lemma} \label{lem:redccwnbhd} Consider the neighborhood of an anticlockwise red loop, with no reddish regions inside. It can be reduced to a diagram where the red loop is absent
(i.e. every region is reddish), or to a diagram as below for $k \ge 3$.

\begin{equation} \label{eq:ccwnbhd}{
\labellist
\small\hair 2pt
 \pinlabel {$v_k$} [ ] at 44 47
\endlabellist
\centering
\ig{1}{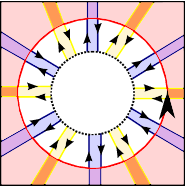}
}\end{equation} \end{lemma}

\begin{proof} Consider the subdiagram which is the interior of the red loop. The boundary of this subdiagram can not contain the color green, lest the exterior of the loop contain a
forbidden brown region. Therefore, the boundary of the interior is either featureless or is a non-empty BS object. Either way, we may apply Claim \ref{BSmaps} to the interior. When the
boundary is featureless, the interior reduces to a box, and then we may apply \eqref{demazureis} or \eqref{bigdemazure} to remove the loop entirely. Otherwise, either the interior is $v_k$
for $k \ge 3$ as desired, or there is a cup/cap on the boundary. In the latter case, we apply \eqref{R2moves} to pull the cup/cap out of the red loop, reducing the size of the
Bott-Samelson boundary by $1$. By induction, this particular anticlockwise loop either disappears or has the desired form. \end{proof}

Note that the diagram \eqref{eq:ccwnbhd} above is an example of a $2$-morphism with an anticlockwise loop which can not be reduced to diagrams without loops. These kinds of diagrams can
appear within the image of $\FC$ applied to non-elliptic webs (c.f. Claim \ref{nonelliptic}). The non-elliptic condition is equivalent to the condition $k \ge 3$.

\subsection{Loop and color removal}
\label{subsec-removal}

We work one color at a time.

\begin{lemma} \label{lem:removered} Any diagram reduces to a diagram without clockwise red loops. Any diagram with featureless reddish boundary reduces to a box. \end{lemma}

\begin{proof} Consider the neighborhood of a red loop. Either no other strands intersect it (i.e. it is featureless), or the strands which intersect it form a non-empty yellow-and-blue BS
object. \igc{1}{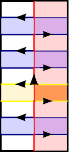} Suppose that this BS object has repetition, i.e. $\mt b \mt b \mt$. Then applying \eqref{circforcesamealt} and \eqref{R2moves}, we can reduce the number of
intersecting strands. \igc{1}{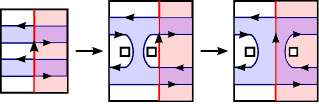} Ergo, up to placing boxes in various regions, we can assume that the intersection is either empty or an alternating BS object. Like every
step of this algorithm, this procedure is local, and will not interfere with other operations.

Consider the case of an innermost clockwise red loop, and let $D$ denote the diagram in its interior. Using the above paragraph, we may assume either that the boundary of $D$ is
featureless, or that it consists of an alternating yellow-and-blue BS object with red added; i.e. it is something like $r p r o r p r o r$ read cyclically around the disk. In the latter
case, we assume that there are $j$ copies of $p$ and $j$ copies of $o$, for some $j>0$. Meanwhile, $D$ has some number $n \ge 0$ of anticlockwise red loops in its interior, so that its
reddish zone appears as below (only red is shown). \igc{1}{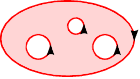} The arguments below will also apply to any diagram $D$ with featureless reddish boundary, regardless of whether or
not it was the interior of a clockwise loop.

We now prove by induction that any such diagram $D$ reduces to a diagram where $n=0$, and that if the diagram has featureless boundary then it reduces to a box. Moreover, if $D$ is the
interior of a red loop and $n=0$, then it reduces to a diagram where $j=0$, i.e. to the featureless case. This will finish the proof of the Lemma: when $n=j=0$ and $D$ is the interior of a
red loop, $D$ reduces to a box which can be slid out using \eqref{slidepolyintoblue} or \eqref{slidepolyintopurple}, and then the loop itself can be eliminated using \eqref{circleis3} or
\eqref{bigcircle}.

We wish to argue that the diagram in the reddish zone can be reduced to a crossingless matching on the $n$-punctured disk (i.e. it has no crossings or closed blue or yellow loops). The
lack of crossings is clear, and any contractible blue or yellow loop can be removed as in the argument of Lemma \ref{lem:allreddish}. Suppose that $n>0$ and there is a blue or yellow
non-contractible loop within the reddish region, and choose an innermost one. \igc{1}{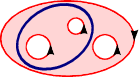} Its interior is another diagram $D'$ with featureless reddish boundary, where the reddish
zone is a crossingless matching. Applying the arguments below to $D'$ will replace $D'$ with a box. Then we can repeat with any other non-contractible loops in a local fashion. Thus we
assume that the reddish zone is a crossingless matching.

Consider the base case when $n=0$. Then the reddish zone is simply-connected. If its boundary is featureless then the reddish zone just contains a box. If not, then any crossingless
matching on a simply-connected region has a cup/cap on the boundary. Using \eqref{R2moves} will lower $j$. To be very precise, this will remove one index from the alternating BS object on
the boundary, producing a repetition, which can then be removed as in the first paragraph; the overall effect is to lower $j$ by one. (This phenomenon will repeat often enough that we
leave it unstated - whenever we remove a single index from the BS boundary, we say that it lowers $j$ by one, leaving unstated the fact that one must then remove the repeated color.) This
handles the base case.

Now suppose that $n>1$. We can assume by Lemma \ref{lem:redccwnbhd} that each anticlockwise red loop has a neighborhood which looks like \eqref{eq:ccwnbhd} for some $k \ge 3$. We call
this a \emph{puncture of size $k$}. We now proceed by induction on the size $k$ of any given puncture, or the number $j$ of the exterior. We do all our computations within the reddish
zone, treating the interiors of the anticlockwise red loops (the punctures) as black boxes.

Consider what any red (not reddish!) region will look like. Instead of reading a list of region colors around a boundary, as we have done in previous arguments, now consider the list of
strand colors around the boundary of the region. They must alternate between red and blue/yellow. If the color blue appears twice (even with yellow in between) applying
\eqref{bigcircforce} will chop the red region into smaller red regions. \igc{1}{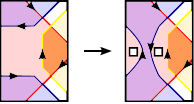} If the red region is a bigon, then one can apply \eqref{R2moves} as usual to pull the cup/cap
either out of the exterior red loop, or into one of the punctures, and use induction. Thus we can assume each red region is a square, having one blue and one yellow strand as parallel
walls. \igc{1}{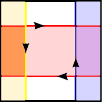} The non-reddish parts of this diagram all occur either inside a puncture or outside the exterior clockwise loop. Thus we can think of red regions as ``edges" which connect punctures or the exterior loop; we will use such a graph below.

Consider what any purple (resp. orange) region will look like. Its boundary alternates between red and blue strands, and it could be a bigon, square, hexagon, octagon, etc. We now argue
that the diagram will simplify unless the purple region has at least 6 sides. If it is a bigon, i.e. if a cup/cap is attached to one of the punctures or the exterior, then applying
\eqref{R2moves} will reduce $k$ or $j$ by one. If it is a square, connecting two punctures or a puncture with the exterior, then one has a local picture like \igc{1}{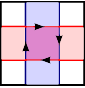} where
the two red strands belong to different red loops. Applying \eqref{squaresbsb}, there are two terms which remain. \igc{1}{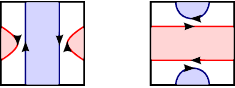} In the second term, the value of $k$ of each
puncture (or $j$ for the exterior) is reduced by one. In the first term, the two punctures are fused into one puncture, or the puncture is fused with the exterior loop, and either case
reduces $n$ by one.

Finally, if the purple region is a square connecting the exterior with itself or connecting any puncture with itself, then one has a local picture like \igc{1}{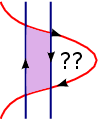} The ???
region can be quite complicated, having punctures, additional strands crossing the red strand, etc. However, once again, applying \eqref{squaresbsb} will yield two terms. In the second
term, the value of either $j$ or some $k$ will decrease. In the first term, either the original exterior loop is split into two smaller loops, or a new clockwise loop is created within a
puncture (containing the ??? region). Either way, the clockwise loop containing the ??? region will be strictly smaller (by one of our inductive criteria) than the original loop, and thus
will disappear. The remainder will also have had either $j$ or some $k$ decreased.

Thus induction applies unless every red region is a square and every purple and orange region is a hexagon or larger. To finish the proof, we make a planar Euler characteristic argument to
prove that this never occurs.

Let $n_k$ be the number of punctures with boundary of size $k$, as above. Let $p_k$ (resp. $o_k$) denote the number of purple (resp. orange) regions with $2k$ sides, and $r$ the number of
red regions. Contract every puncture into a vertex, and if $j>0$ then contract the outer loop into a vertex, so that one obtains a graph embedded in $S^2$. Now, each puncture of size $k$
yields a $4k$-valent vertex, with $2k$ red regions, $k$ purple regions, and $k$ orange regions adjoining. After this contraction, the purple regions enumerating $p_k$ have $k$ adjoining
edges, and the red regions are bigons. Therefore, by counting purple (resp. orange) regions paired with an adjoining vertex, we obtain \begin{equation} \label{euler1} \sum_k k p_k = j +
\sum_k k n_k = \sum_k o_k. \end{equation} Counting red regions paired with a vertex, we obtain \begin{equation} \label{euler2} 2r = \sum_k 2k n_k + 2j, \end{equation} so that $r$ is also
equal to the quantity in \eqref{euler1}.

Each edge appears adjacent to a single red region, so the number of edges is $2r$. The number of vertices is $1+\sum n_k$ if $j > 0$, and $\sum n_k$ if $j=0$. The number of regions is $r +
\sum p_k + \sum o_k$. Thus (when $j>0$) we have: \[ 2 = V - E + R = (1-j) + \sum (1-k)n_k + \sum p_k + \sum o_k. \] When $j=0$, the $(1-j)$ term is ignored. Now, by assumption each $k \ge
3$. Therefore \[ \sum p_k + \sum o_k \le \frac{1}{3} \sum (kp_k + k o_k) = \frac{2}{3} (j+ \sum k n_k), \] so that \[2 \le (1 - \frac{j}{3}) + \sum (1 - \frac{k}{3}) n_k.\] However,
$(1-\frac{k}{3})$ is nonpositive, and the term $(1-\frac{j}{3})$ is at most one (and disappears when $j=0$). This is a contradiction, and concludes the proof of Lemma \ref{lem:removered}.
\end{proof}

\begin{cor} Let $K \subsetneq S$. Then $\iota_K$ is full. \end{cor}

\begin{proof} Suppose that $K$ does not contain the primary color red. Then any $2$-morphism between $1$-morphisms in $\DC(K)$ will have no reddish exterior regions. By removing all
clockwise red loops, we thereby remove the primary color red entirely. The resulting $2$-morphism is in the image of $\iota_K$. \end{proof}

The lemma essentially gave an algorithm to remove any given innermost clockwise loop from a diagram. However, the procedure did screw up the topology of the differently-colored strands
which intersected that loop, and thus can create (blue or yellow) loops as well as destroy them. There does not seem to be an easy way to deduce that one can remove all clockwise loops
of different colors at once, directly from the lemma. Instead, we use a circuitous route.

\subsection{Arguments from categorification}
\label{subsec-catfnargs}

\begin{lemma} \label{niceform} $2$-morphisms in $\DC$ between BS objects are spanned by diagrams where the only instances of purple (resp. orange, green) have neighborhoods which look like
$v_k$ for $k \ge 3$. \end{lemma}

\begin{proof} Consider a purple region in such a diagram, which must be internal, and its interior $D$. Then $D$ is a diagram with featureless reddish boundary, so by Lemma
\ref{lem:removered} it reduces to a box. In this fashion, we can assume every purple region is simply connected, and has no yellowish regions in its neighborhood. This neighborhood is now
a morphism in $\DC(r,b)$, and the result follows from the Cathedral. \end{proof}

Recall that each $v_k$ for $k > 3$ can be rewritten using $v_3$.

BS objects form a monoidal category, which is supposed to categorify the Hecke algebra $\HB$ of $W$. In \cite{EWGR4SB}, a diagrammatic category $\DC_{BS}$ (called $\DC$ in that
paper) is described by generators and relations (for an arbitrary Coxeter group), which is meant to encode the morphisms between BS bimodules. In particular, there is a $2$-functor $\IC$
from $\DC_{BS}$ to $\DC$, defined on dihedral parabolic subgroups in \cite{EDihedral}, whose $1$-morphism image consists of the BS objects. The $2$-morphism image consists of those
morphisms between BS objects generated by cups and caps, as well as the map $v_3$.

\begin{prop} \label{prop:bs-to-sbs} The natural $2$-functor $\IC$ from $\DC_{BS}$ to $\DC$ is fully faithful. The $2$-functor $\GC$ is faithful from objects in the image of $\IC$.
\end{prop}

\begin{proof} Lemma \ref{niceform} indicates that $\IC$ is full. The composition $\GC \circ \IC$ was shown to be faithful in \cite{EWGR4SB}, so that $\IC$ must also be faithful, and
$\GC$ must be faithful from the image. \end{proof}

The wonderful fact about Soergel bimodules and their diagrammatic categories is that the behavior of morphism spaces is governed by the BS objects.

\begin{cor} \label{cor:maincor} The $2$-functor $\GC$ is faithful. $2$-morphism spaces in $\DC$ satisfy the Soergel-Williamson Hom Formula \cite[Theorem 7.2.2]{WilSSB}, and $\DC$
categorifies the Hecke algebroid. \end{cor}

\begin{proof} This is a quick summary of an argument made almost verbatim in the Cathedral. Recall the definitions of the \emph{Hecke algebra} $\HB$ (resp. \emph{Hecke algebroid} $\HG$)
from \cite{WilSSB} or from \cite[section 2.4]{EDihedral}. Recall the definition of a \emph{potential categorification} of the Hecke algebra (resp. algebroid) from \cite[sections 2.3 and
2.4.3]{EDihedral}. Proposition \ref{prop:2colorGC}) implies that $\DC$ categorifies all the relations in $\HG$ coming from dihedral parabolic subgroups. As all the relations in $\HG$ come
from finite parabolic subgroups, we see that $\DC$ is a potential categorification of $\HG$, and therefore induces a trace map on $\HG$. Any trace map on $\HG$ is determined by its
evaluation on the subalgebroid $\HB$, or in other words, the sizes of arbitrary morphism spaces in $\DC$ (and their behavior under $\GC$) are determined by the morphism spaces between BS
objects. Thus Proposition \ref{prop:bs-to-sbs} implies that $\GC$ is faithful everywhere, as desired (c.f. \cite[sections 2.4.3 and 5.4.2]{EDihedral}). \end{proof}

Theorem \ref{superthm-sl3} follows immediately. So does the fact that all clockwise loops can be removed simultaneously. One may also assume (up to reduction) that every region
in a diagram is simply-connected.

\subsection{Proof of Theorem \ref{thm:sl3main}}
\label{subsec-thmpf}

Recall the definition of $\mSBSBim$ and the functor $\FC \co \FundOm \to \mSBSBim \subset \DC$ from \S\ref{subsec-equiv-sl3}. Our goal is to show that $\FC$ is a degree zero equivalence.
Much of what needs to be shown is already implied by Theorem \ref{superthm-sl3}. It is an exercise in the Hecke algebroid and the basic properties of $\SSBim$ to show that $\FC$ is
essentially surjective up to grading shift.

We wish to show that $\FC$ is fully faithful in degree zero. However, morphism spaces in $\DC$ have a graded dimension governed by the Soergel-Williamson Hom Formula. The fact that
morphisms between $1$-morphisms in the image of $\FC$ are positively graded, and that their degree $0$ morphisms have the correct size, is a purely numerological one. It was first observed
by Lusztig in the seminal paper \cite{LusztigGS} from which the idea of Geometric Satake was birthed. Thus it remains to show that $\FC$ is full, which is a purely diagrammatic argument.

A diagram representing a morphism in $\mSBSBim$ has no exterior white regions. Consider a given white region. Using Proposition \ref{prop:colorremoval} we can remove any color on the
interior of the white region. This leaves a box, and nothing more. The region is delineated by $k$ primary-colored strands, with no two of the same color being adjacent. If $k=1$ the
region can be removed using \eqref{demazureis}. If $k=2$ the region can be removed using \eqref{singR2nonoriented1}. If $k=3$ then a neighborhood of this triangular white region is
precisely $\FC$ applied to the trivalent vertex \igc{1}{3FundtoSingSoerg} though possibly with a box in the white region. If $k>4$ then there must be a repetition of colors. Using a
now-familiar argument, one applies \eqref{circforcesamealt} to this repetition, which splits the white region into two smaller regions. Therefore, we may assume that every white region is
a triangle.

Remove a neighborhood of each white region, and consider the remaining diagram on the punctured disk. The remaining diagram has no crossings, lest there exist a brown region. We wish to
show that it is in the image of $\FC$ applied to a diagram on the punctured disk without any trivalent vertices. In other words, we must show that every red region has exactly one purple
and one orange region adjacent to it (and something similar for blue and yellow regions). Note that, because of the behavior of the diagram at the boundary and near each puncture, the
regions surrounding the red region must alternate between orange and purple.

Consider any red region. Exactly as in the proof of Lemma \ref{lem:removered}, if it abuts two purple regions then one may fuse them with \eqref{bigcircforce}, splitting the red region
into two smaller regions. This can be repeated until each red region has exactly one purple and one orange region adjacent to it, as desired.

Therefore, every morphism between objects in $\mSBSBim$ is in the span of diagrams in the image of $\FC$, but with boxes in various regions. Only when there are no boxes will the diagram
have degree zero. This concludes the proof of fullness. In fact, we have proven a slightly stronger statement.

\begin{prop} Morphism spaces (of arbitrary degree) in $\mSBSBim$ between objects in the image of $\FC$ are generated by diagrams in the image of $\FC$, under the operation of placing boxes
in various regions. \end{prop}

\bibliographystyle{plain}
\bibliography{mastercopy}{}

\end{document}